\documentclass[a4paper,10pt,reqno]{amsart}

\usepackage[top=2.5cm, bottom=2cm, left=2.5cm, right=2.5cm]{geometry}

\usepackage{amssymb,amsmath,amsthm,ascmac,array,multirow}
\usepackage[dvipdfmx]{graphicx}
\usepackage{color}

\allowdisplaybreaks

\newtheorem{thm}{Theorem}[section]

\newtheorem{lem}[thm]{Lemma}

\newtheorem{rem}[thm]{Remark}
\newtheorem{ass}[thm]{Assumption}

\newcommand{\mcf}{\mathcal{F}}

\newcommand{\mcl}{\mathcal{L}}

\newcommand{\mfp}{\mathfrak{p}}

\newcommand{\mfM}{\mathfrak{M}}

\newcommand{\mbbh}{\mathbb{H}}
\newcommand{\mbbn}{\mathbb{N}}
\newcommand{\mbbr}{\mathbb{R}}
\newcommand{\mbbu}{\mathbb{U}}

\newcommand{\mbbz}{\mathbb{Z}}
\newcommand{\mbX}{\mathbf{X}}
\newcommand{\mbY}{\mathbf{Y}}

\newcommand{\mbx}{\mathbf{x}}
\newcommand{\mby}{\mathbf{y}}

\newcommand{\D}{\Delta}
\newcommand{\Sig}{\Sigma}

\newcommand{\Gam}{\Gamma}
\newcommand{\p}{\partial}

\newcommand{\cil}{\xrightarrow{\mcl}} 

\newcommand{\cip}{\xrightarrow{P}} 

\newcommand{\argmin}{\mathop{\rm argmin}} 
\newcommand{\argmax}{\mathop{\rm argmax}}

\def\sumj{\sum_{j=1}^{n}}

\def\cov{\mathrm{cov}}
\def\dim{\mathrm{dim}}

\def\tz{\theta_{0}}

\def\tes{\hat{\theta}_{n}}

\def\qbic{\mathrm{QBIC}}

\title{Model comparison for generalized linear models with dependent observations}

\author{Shoichi Eguchi}
\address{Graduate School of Mathematics, Kyushu University. 
744 Motooka, Nishi-ku, Fukuoka 819-0395, Japan.}
\email{s-eguchi@math.kyushu-u.ac.jp}

\date{\today}

\begin{document}
\setlength{\baselineskip}{5mm}

\maketitle

\begin{abstract}
The stochastic expansion of the marginal quasi-likelihood function associated with a class of generalized linear models is shown.
Based on the expansion, a quasi-Bayesian information criterion is proposed that is able to deal with misspecified models and dependent data, resulting in a theoretical extension of the classical Schwarz's Bayesian information criterion.
It is also proved that the proposed criterion has model selection consistency with respect to the optimal model.
Some illustrative numerical examples and a real data example are presented.
\end{abstract}

\section{Introduction} \label{Intro}


The generalized linear model (GLM, McCullagh and Nelder  \cite{McCNel83}) is an extension of the linear regression model and has many applications and extensions;
for example, actuarial science (Antonio and Beirlant \cite{AntBei07}, Haberman and Renshaw \cite{HabRen96}), GLMixedM in risk management (McNeil and Wendin \cite{McNWen07}), and generalized additive models (Berg \cite{Ber07}, Hastie and Tibshirani \cite{HasTib90}).

We consider data $(y_{j},x_{j})_{j=1}^{n}=(y_{j},x_{j,1},\ldots,x_{j,p})_{j=1}^{n}$, where the $y_{j}$'s and $x_{j}$'s are realizations of the response variables $\mbY_{n}=(Y_{1},\ldots,Y_{n})^{\prime}$ and the explanatory variables $\mbX_{n}=(X_{1},\ldots,X_{n})^{\prime}$, respectively; here, the notation $\prime$ means the transpose.  
Furthermore, we will assume that the conditional distribution of $\mbY_{n}$ given $\mbX_{n}$ is given by a GLM.
Then, the conditional distribution is assumed to belong to an exponential family, for example, normal, binomial, or Poisson distributions.
In this paper, a result is presented about the stochastic expansion of the marginal quasi-likelihood function associated with a class of possibly misspecified GLMs for dependent data. 
Based on this expansion, we propose the quasi-Bayesian information criterion, which is an extension of the generalized BIC given by Luv and Liu \cite{LvLiu14}.

The Bayesian principle for model selection is formulated as follows.
Suppose that $M$ Bayesian candidate models $\mfM_{1},\ldots,\mfM_{M}$ are given. 
Each $\mfM_{m}$ is described by $\big\{\big(\mfp_{m},\pi_{m}(\theta),\mbbh_{m,n}(\theta)\big)\big|\theta\in\Theta_{m}\big\}$, where $\mfp_{m}$ is the non-zero prior relative occurrence probability of the $m$th model of the $M$ models, $\pi_{m}$ is the prior-probability density on $\Theta_{m}$, and $\mbbh_{m,n}$ is the logarithmic quasi-likelihood function.
Here ``quasi" means that a parametric model that may not include a true distribution of the observation data is considered.
The conventional Bayesian principle of model selection for $\mfM_{1},\dots,\mfM_{M}$ chooses the model that is most likely in terms of the posterior probability, i.e., choosing the model that maximizes $P(\mfM_{m}|\mby_{n})$, where
\begin{align*}
P(\mfM_{m}|\mby_{n})&=\frac{\Big(\int_{\Theta_{m}}\exp\{\mbbh_{m,n}(\theta)\}\pi_{m}(\theta)d\theta\Big)\mfp_{m}}{\sum_{i=1}^{M}\Big(\int_{\Theta_{i}}\exp\{\mbbh_{i,n}(\theta)\}\pi_{i}(\theta)d\theta\Big)\mfp_{i}},
\end{align*}
where $\int_{\Theta_{m}}\exp\{\mbbh_{m,n}(\theta)\}\pi_{m}(\theta)d\theta$ is called the marginal quasi-likelihood function.
When the prior plausibilities on the $M$ competing models are equal, we select the model that maximizes the marginal quasi-likelihood function; even if the prior probabilities are not equal, we can trivially correct the selection manner by the factors $\mfp_{m}$. 
Hence, the logarithm of the marginal quasi-likelihood function 
\begin{align*}
\log\bigg(\int_{\Theta_{m}}\exp\{\mbbh_{m,n}(\theta)\}\pi_{m}(\theta)d\theta\bigg)
\end{align*}
is used to select the model.
 
As explained in \cite{LvLiu14}, another interpretation of model selection is possible through the Kullback--Leibler divergence (KL divergence). 
The KL divergence between true conditional model $g_{n}$ and marginal quasi-likelihood function $\int_{\Theta_{m}}$ $\exp\{\mbbh_{m,n}(\theta)\}\pi_{m}(\theta)d\theta$ is given by
\begin{align}
&I\bigg(g_{n};\int_{\Theta_{m}}\exp\{\mbbh_{m,n}(\theta)\}\pi_{m}(\theta)d\theta)\bigg) \notag \\
&=E[\log g_{n}(\mbY_{n}|\mbX_{n})]+E\bigg[-\log\bigg(\int_{\Theta_{m}}\exp\{\mbbh_{m,n}(\theta)\}\pi_{m}(\theta)d\theta\bigg)\bigg], \label{UnbiEst}
\end{align}
where the expectation is taken with respect to true distribution $G_{n}$. 
Because of (\ref{UnbiEst}), $-\log\big(\int_{\Theta_{m}}\exp$ $\{\mbbh_{m,n}(\theta)\}\pi_{m}(\theta)d\theta\big)$ is clearly an unbiased estimator of $I\big(g_{n};\int_{\Theta_{m}}\exp\{\mbbh_{n}(\theta;\cdot)\}\pi(\theta)d\theta\big)$ except for a constant term free of $\theta$.
Hence, the Bayesian principle of model selection can be restated as choosing the model that minimizes the KL divergence of the marginal quasi-likelihood function from the true distribution.
Note that (\ref{UnbiEst}) holds regardless of whether or not the set of candidate models contains the true model.
In particular, assume that $\mbX_{n}$ is absent and that $\mbbh_{m,n}(\theta)=\sumj\log $ $f_{m,n}(y_{j};\theta)$ for the case of independent observations $y_{1},\ldots,y_{n}$ with correctly specified exponential family models. 
Then, Schwarz \cite{Sch78} showed that the logarithmic marginal quasi-likelihood $\log(\int_{\Theta_{m}}\exp\{\mbbh_{m,n}(\theta)\}\pi_{m}(\theta)d\theta)$ admits the stochastic expansion
\begin{align}
\log\bigg(\int_{\Theta_{m}}\exp\{\mbbh_{m,n}(\theta)\}\pi_{m}(\theta)d\theta\bigg)&=\sumj\log f_{m,n}(y_{j};\hat{\theta}_{m,n}^{\mathrm{MLE}}) \notag \\
&\qquad-\frac{p}{2}\log n+O_{p}(1), \label{shw}
\end{align}
with $\hat{\theta}_{m,n}^{\mathrm{MLE}}$ denoting the maximum likelihood estimator of $\theta$, under some regularity conditions. 
Because of (\ref{shw}), we obtain the classical Bayesian information criterion for model selection:
\begin{align*}
\mathrm{BIC}=-2\sumj\log f_{m,n}(y_{j};\hat{\theta}_{m,n}^{\mathrm{MLE}})+p\log n. 
\end{align*}

In the past, many authors have investigated several kinds of information criteria for model selection in various settings; see, for example, Burnham and Anderson \cite{BurAnd02} for an account of these developments.
Bozdogan \cite{Boz87} showed that the Akaike information criterion (AIC, Akaike \cite{Aka73}, \cite{Aka74}) has a positive probability of overestimating the true dimension.
Casella {\it et al.} \cite{CGMM09} and Fasen and Kimmig \cite{FaKim15} as well as the references therein studied the model selection consistency of the BIC.
Moreover, various extensions of the AIC and BIC have been introduced; for example, the extended BIC for large model spaces (Chen and Chen \cite{CheChe08}), generalized information criterion (Konishi and Kitagawa \cite{KoniKita96}), generalized BIC in misspecified GLMs for independent data (Lv and Liu \cite{LvLiu14}), and information criteria for stochastic processes (e.g., Sei and Komaki \cite{SeiKom07} and Uchida \cite{Uch10}). 

The rest of the paper is organized as follows. In Section \ref{NoTe}, we describe our working model, notations, and assumptions. 
We also discuss the asymptotic properties of the quasi-maximum likelihood estimator in possibly misspecified GLMs for dependent data.
Section \ref{qbic.sec} presents the stochastic expansion of the logarithmic marginal quasi-likelihood and model selection consistency with respect to the optimal model (see Section 3.2 for the definition). 
In Section \ref{ExSim}, we illustrate the performance of the model selection criteria in correctly specified and misspecified models.
Section \ref{real.simu} presents a real data example. 
The proofs of our results are given in Section \ref{Proofs} and Supplementary Material.


\section{Quasi-maximum likelihood estimation of dependent GLM} \label{NoTe}

Let $\mbY_{n}=(Y_{1},\ldots,Y_{n})^{\prime}$ be the $n$-dimensional random vector and $\mbX_{n}=(X_{1},\ldots,X_{n})^{\prime}$ be the $n\times p$ random time series. 
We write $X_{j}=(X_{j,1},\ldots,X_{j,p})^{\prime}$. 
We assume that the unknown true distribution of $(\mbX_{n},\mbY_{n})$ has the density $g_{n}$ with respect to some dominating $\sigma$-finite measure:
\begin{align*}
g_{n}(\mbx_{n},\mby_{n})=g_{n}(\mbx_{n})g_{n}(\mby_{n}|\mbx_{n}),
\end{align*}
where $\mbx_{n}=(x_{1},\ldots,x_{n})^{\prime}$, $x_{j}=(x_{j,1},\ldots,x_{j,p})^{\prime}$ and $\mby_{n}=(y_{1},\ldots,y_{n})^{\prime}$.


\subsection{Model setup}
We consider possibly misspecified $M$ candidate models to estimate the true model $G_{n}$. 
Assume that each candidate model is given by
\begin{align}
f_{m,n}(\mbx_{n},\mby_{n};\theta)=f_{n}(\mbx_{n})f_{m,n}(\mby_{n}|\mbx_{n};\theta)=f_{n}(\mbx_{n})\prod_{j=1}^{n}f_{m,n,j}(y_{j}|x_{j};\theta),
\label{candi_join}
\end{align}
with $\theta=(\theta_{1},\ldots,\theta_{p_{m}})\in\Theta_{m}$, where the $m$th parameter space $\Theta_{m}\subset\mbbr^{p_{m}}$ is a bounded convex domain and $p_{m}\leq p$. 
Although the true structure is unknown, (\ref{candi_join}) means that the candidate models are given as if $Y_{1},\ldots,Y_{n}$ are $(X_{1},\ldots,X_{n})$-conditionally independent and each $(X_{1},\ldots,X_{n})$-conditional distribution of $Y_{j}$ only depends on $X_{j}$.
The later condition is used to simplify the theoretical consideration and is not essential.
Because (\ref{candi_join}) entails that the candidate distribution of $\mbX_{n}$ does not depend on the parameter, only the conditional distribution of $\mbY_{n}$ given $\mbX_{n}$ is considered.
GLM $\mfM_{m}$ is used as the working model, with respect to some dominating measure: $\mfM_{m}$ is given by
\begin{align}
f_{m,n}(\mby_{n}|\mbx_{n};\theta)&=\prod_{j=1}^{n}f_{m,n,j}(y_{j}|x_{j};\theta) \notag \\
&=\prod_{j=1}^{n}\exp\big(y_{j}x_{j}^{\prime}\theta-b_{m}(x_{j}^{\prime}\theta)+c_{m}(y_{j})\big),
\label{glm_pdf}
\end{align}
where, for brevity, we write $x_{j}^{\prime}\theta=\sum_{i=1}^{p_{m}}x_{j,d_{i}(m)}\theta_{i}$ with $\{d_{1}(m),\ldots,d_{p_{m}}(m)\}\subset\{1,\ldots,p\}$ for any $m$.
Further, $b_{m}(\cdot)$ and $c_{m}(\cdot)$ are determined by each assumed conditional distribution of $\mbY_{n}$ given $\mbX_{n}$, and $b_{m}(\cdot)$ is a sufficiently smooth convex function defined on $\mbbr$. For example, $b_{m}(\theta)=\theta^{2}/{2}$ in the case of Gaussian regression and $b_{m}(\theta)=\log(1+e^{\theta})$ in the case of logistic regression. 
We assume that $b_{1}(\theta)=\dots=b_{M}(\theta)$ and $c_{1}(y)=\dots=c_{M}(y)$ for simplicity, that is, only consider variable selection concerning $\mbX_{n}$.

Because each candidate model $\mfM_{m}$ is possibly misspecified and $c_{m}(\cdot)$ of (\ref{glm_pdf}) is independent of $\theta$, the logarithmic quasi-likelihood function $\mbbh_{m,n}$ can be defined by
\begin{align}
\mbbh_{m,n}(\theta)=\sumj\big(Y_{j}X_{j}^{\prime}\theta-b_{m}(X_{j}^{\prime}\theta)\big).
\label{Qlf}
\end{align}
The quasi-maximum likelihood estimator (QMLE) associated with $\mbbh_{m,n}$ is defined to be any maximizer of $\mbbh_{m,n}$:
\begin{align*}
\hat{\theta}_{m,n}\in\argmax_{\theta\in\Theta_{m}}\mbbh_{m,n}(\theta).
\end{align*}
Clearly, when $b_{m}$ is differentiable, $\hat{\theta}_{m,n}$ is the solution to the quasi-score function
\begin{align*}
\p_{\theta}\mbbh_{m,n}(\theta)=\sumj\big(Y_{j}-\p b_{m}(X_{j}^{\prime}\theta)\big)X_{j}=0,
\end{align*}
where $\p_{\theta}=\p/\p\theta$ and $\p b_{m}(x^{\prime}\theta)=\frac{\p}{\p(x^{\prime}\theta)}b_{m}(x^{\prime}\theta)$. 

From now on, we will omit the model index $``m"$ from the notation for notational brevity. 


\subsection{Asymptotic behavior of the QMLE}
In this section, we will show the asymptotic properties of the QMLE. 
Fahrmeir and Kaufmann \cite{FahKau} studied the consistency and asymptotic normality of the MLE in correctly specified GLMs.
Moreover, White \cite{Whi82} and Domowitz and White \cite{DomWhi82} investigated the properties of the QMLE in misspecified models and treated independent data and dependent data, respectively.
The settings of Domowitz and White \cite{DomWhi82} are more general than the currently considered settings.
In contrast, we will assume clearer conditions that are more tailored for the GLMs. 

Denote by $\mcf_{j}=\sigma(Y_{i},X_{i};i\leq j)$ the $\sigma$-field representing the data information at stage $j$.
If $a_{n}$ and $b_{n}$ satisfy $a_{n}\leq Cb_{n}$ for some constant $C>0$, we write $a_{n}\lesssim b_{n}$.
We assume the following conditions.

\begin{ass}
For some constant $C\geq0$ and $C^{\prime}\geq0$, 

(i) $\displaystyle \max_{i\in\{1,2,3\}}|\p^{i}b(x)|\lesssim 1+|x|^{C}$,

(ii) $\displaystyle E\big[|Y_{j}|^{3}|\mcf_{j-1}\vee\sigma(X_{j})\big]\lesssim 1+|X_{j}|^{C^{\prime}}$ $a.s.$ for any $j\in\mbbn$,

(iii) $\displaystyle \sup_{j\in\mbbn}E\big[|X_{j}|^{3C+C^{\prime}+3}\big]<\infty$.
\label{Ass1}
\end{ass}

\begin{ass}
There exists a measurable function $F:\mbbr^{p}\to\mbbr$ such that $E[Y_{j}|\mcf_{j-1}\vee\sigma(X_{j})]=F(X_{j})$ for every $j\in\mbbn$.
\label{Ass2}
\end{ass}

\begin{ass}
Let $\zeta_{j}=(X_{j},Y_{j})$ for any $j$. For some $c>0$, 
\begin{align*}
\alpha(k)\leq c^{-1}e^{-ck}
\end{align*}
for all $k\in\mbbn$, where
\begin{align*}
\alpha(k):=\sup_{j\in\mbbn}\sup_{\substack{A\in\sigma(\zeta_{i};i\leq j) \\ B\in\sigma(\zeta_{i};i\geq j+k)}}|P[A\cap B]-P[A]P[B]|.
\end{align*}
\label{Ass5}
\end{ass}

When Assumption \ref{Ass5} holds, $\{\zeta_{j};j=1,2,\ldots\}$ is called exponential $\alpha$-mixing.
In particular, Assumption \ref{Ass5} implies that $\psi_{j}:=\big(Y_{j}-F(X_{j})\big)X_{j}$, $j\in\mbbn$, is exponential $\alpha$-mixing.

\begin{ass}
There exists a non-degenerate probability measure $\nu$ such that the following holds:
\smallskip

(i) $\displaystyle\frac{1}{n}\sumj\big(F(X_{j})X_{j}^{\prime}\theta-b(X_{j}^{\prime}\theta)\big)-\int\big(F(x)x^{\prime}\theta-b(x^{\prime}\theta)\big)\nu(dx)\cip0$ for any $\theta\in\Theta$, 
\smallskip

(ii) $\displaystyle\frac{1}{n}\sumj\p^{2}b(X_{j}^{\prime}\theta)X_{j}X_{j}^{\prime}-\int\p^{2}b(x^{\prime}\theta)xx^{\prime}\nu(dx)\cip0$ for any $\theta\in\Theta$, \\
where the notation $\cip$ means the convergence in probability.
\label{Ass4}
\end{ass}

\begin{ass}
There exists a function $\underline{b}:\mbbr^{p}\to(0,\infty)$, (i) for any $x$, $\displaystyle\inf_{\theta\in\Theta}\p^{2}b(x^{\prime}\theta)\geq\underline{b}(x)$,

(ii) for some constant $\lambda_{0}>0$, $\displaystyle\limsup_{n\to\infty}P\bigg[\lambda_{\min}\bigg(\frac{1}{n}\sumj\underline{b}(X_{j})X_{j}X_{j}^{\prime}\bigg)<\lambda_{0}\bigg]=0$, where $\lambda_{\min}(\cdot)$ denotes the smallest eigenvalue of a given matrix. 
\label{Ass3}
\end{ass}

These assumptions are explained below.
Assumptions \ref{Ass1}, \ref{Ass2}, and \ref{Ass5} ensure the asymptotic properties (see Lemmas \ref{lem3} and \ref{lem1}).
Moreover, these assumptions facilitate the derivation of the consistency and asymptotic normality of the QMLE $\tes$, which are given in Theorems \ref{consis} and \ref{AN}, respectively.

If Assumptions \ref{Ass1} and \ref{Ass2} hold, for some constants $C\geq0$ and $C^{\prime}\geq0$,
\begin{align}
&\sup_{n>0}E\left[\sup_{\theta\in\Theta}\left|\p_{\theta}\left\{\frac{1}{n}\sumj\big(F(X_{j})X_{j}^{\prime}\theta-b(X_{j}^{\prime}\theta)\big)\right\}\right|\right] \notag \\
&\leq\sup_{n>0}\frac{1}{n}\sumj E\bigg[\Big(\big|F(X_{j})\big|+\sup_{\theta\in\Theta}\big|b(X_{j}^{\prime}\theta)\big|\Big)|X_{j}|\bigg] \notag \\
&\leq\sup_{n>0}\frac{1}{n}\sumj E\bigg[\Big(E[|Y_{j}||\mcf_{j-1}\vee\sigma(X_{j})]+(1+|X_{j}|)^{C}\Big)|X_{j}|\bigg] \notag \\
&\lesssim\sup_{n>0}\frac{1}{n}\sumj E\bigg[\Big((1+|X_{j}|)^{C^{\prime}}+(1+|X_{j}|)^{C}\Big)|X_{j}|\bigg]<\infty.
\label{remeq2}
\end{align}
Because (\ref{remeq2}) gives the tightness of $\big\{\frac{1}{n}\sumj\big(F(X_{j})X_{j}^{\prime}\theta-b(X_{j}^{\prime}\theta)\big)\big\}$ in $\mathcal{C}(\Theta)$, we have
\begin{align}
\sup_{\theta\in\Theta}\left|\frac{1}{n}\sumj\big(F(X_{j})X_{j}^{\prime}\theta-b(X_{j}^{\prime}\theta)\big)-\int\big(F(x)x^{\prime}\theta-b(x^{\prime}\theta)\big)\nu(dx)\right|\cip0
\label{unicon1}
\end{align}
under Assumptions \ref{Ass1}, \ref{Ass2}, and \ref{Ass4} (i).

From Assumptions \ref{Ass1}--\ref{Ass5} and \ref{Ass4} (i), we have
\begin{align}
\frac{1}{n}\mbbh_{n}(\theta)&=\frac{1}{n}\sumj\psi_{j}^{\prime}\theta+\frac{1}{n}\sumj\big\{\big(F(X_{j})X_{j}^{\prime}\theta-b(X_{j}^{\prime}\theta)\big)\big\} \notag \\
&=O_{p}\bigg(\frac{1}{\sqrt{n}}\bigg)+\frac{1}{n}\sumj\big\{\big(F(X_{j})X_{j}^{\prime}\theta-b(X_{j}^{\prime}\theta)\big)\big\} \notag \\
&\cip\int\big(F(x)x^{\prime}\theta-b(x^{\prime}\theta)\big)\nu(dx)=:\mbbh_{0}(\theta). \label{lawL}
\end{align}
The proof of the tightness of $\big\{\frac{1}{\sqrt{n}}\sumj\psi_{j}\big\}$ used above is given in Lemma \ref{lem3}.

Assumption \ref{Ass3} ensures the uniqueness of the QMLE $\tes$ by the strict concavity of $\mbbh_{n}$.
Because Assumption \ref{Ass4} (ii) gives $-\frac{1}{n}\p_{\theta}^{2}\mbbh_{n}(\theta)=-\p_{\theta}^{2}\mbbh_{0}(\theta)+o_{p}(1)$ for any $\theta$, Assumptions \ref{Ass1}--\ref{Ass3} imply that the equation 
\begin{align*}
\p_{\theta}\mbbh_{0}(\theta)=\int\big(F(x)-\p b(x^{\prime}\theta)\big)x\nu(dx)=0
\end{align*} 
admits a unique solution.
The {\it optimal} parameter $\tz$ may be defined as the unique maximizer of $\mbbh_{0}(\theta)$:
\begin{align*}
\{\tz\}=\argmax_{\theta\in\Theta}\mbbh_{0}(\theta).
\end{align*}
The quasi-observed information is given by $\Gam_{n}:=-\frac{1}{n}\p_{\theta}^{2}\mbbh_{n}(\tz)=\frac{1}{n}\sumj$ $\p^{2}b(X_{j}^{\prime}\tz) X_{j}X_{j}^{\prime}$ so that $\Gam_{n}$ satisfies the equation 
\begin{align*}
\Gam_{n}=\Gam_{0}+o_{p}(1),
\end{align*}
where $\Gam_{0}:=\int\p^{2}b(x^{\prime}\tz)xx^{\prime}\nu(dx)$.

\begin{rem}{\rm
The $\beta$-mixing coefficients of $\{\zeta_{j}\}$ are defined by
\begin{align*}
\beta(k):=\sup_{j\in\mbbn}E\bigg[\sup_{B\in\sigma(\zeta_{i};i\geq j+k)}\big|P\big(B|\sigma(\zeta_{i};i\leq j)\big)-P(B)\big|\bigg].
\end{align*} 
If $\beta(k)=O(e^{-ak})$ for some $a>0$ and for all $k\in\mbbn$, then $\{\zeta_{j}\}$ is called exponential $\beta$-mixing (e.g., Davydov \cite{Dav73} and Liebscher \cite{Lie05}).
The exponential $\beta$-mixing property implies the exponential $\alpha$-mixing property.
When we replace Assumption \ref{Ass5} with the exponential $\beta$-mixing property of $\{\zeta_{j}\}$ under some appropriate moment condition, the following conditions follow when  an obvious discrete-time counterpart of Masuda \cite[Lemma 4.3]{Mas13} is applied: 
(i) For some constant $\beta_{1}>0$ and $q_{1}>0$,
\begin{align*}
&\sup_{n>0}E\bigg[\bigg(n^{\beta_{1}}\sup_{\theta\in\Theta}\bigg|\frac{1}{n}\sumj\big(F(X_{j})X_{j}^{\prime}\theta-b(X_{j}^{\prime}\theta)\big) \\
&\qquad\qquad\qquad\qquad-\int\big(F(x)x^{\prime}\theta-b(x^{\prime}\theta)\big)\nu(dx)\bigg|\bigg)^{q_{1}}\bigg]<\infty.
\end{align*}
(ii) For some constant $\beta_{2}>0$ and $q_{2}>0$,
\begin{align*}
\sup_{n>0}E\bigg[\bigg(n^{\beta_{2}}\sup_{\theta\in\Theta}\bigg|\frac{1}{n}\sumj\p^{2}b(X_{j}^{\prime}\theta)X_{j}X_{j}^{\prime}-\int\p^{2}b(x^{\prime}\theta)xx^{\prime}\nu(dx)\bigg|\bigg)^{q_{2}}\bigg]<\infty.
\end{align*}
Because of the Borel-Canntelli lemma, if $q_{1}$ and $q_{2}$ can be taken large enough, we may deduce that almost surely
\begin{align*}
\sup_{\theta\in\Theta}\bigg|\frac{1}{n}\sumj\big(F(X_{j})X_{j}^{\prime}\theta-b(X_{j}^{\prime}\theta)\big)-\int\big(F(x)x^{\prime}\theta-b(x^{\prime}\theta)\big)\nu(dx)\bigg|\to0
\end{align*}
and
\begin{align*}
\sup_{\theta\in\Theta}\bigg|\frac{1}{n}\sumj\p^{2}b(X_{j}^{\prime}\theta)X_{j}X_{j}^{\prime}-\int\p^{2}b(x^{\prime}\theta)xx^{\prime}\nu(dx)\bigg|\to0.
\end{align*}
}\qed
\end{rem}

\begin{rem}{\rm
We can relax Assumption \ref{Ass2} by replacing $E[Y_{j}|\mcf_{j-1}\vee\sigma(X_{j})]=F(X_{j})$ with $E[Y_{j}|\mcf_{j-1}\vee\sigma(X_{j})]=F(X_{j-m+1},\ldots,X_{j})$ for some $m\geq1$.
Then, Assumption \ref{Ass4} (i) is modified.
For example, in the case of $m=2$, Assumption \ref{Ass4} (i) could be replaced by
\begin{align*}
&\frac{1}{n}\sum_{j=2}^{n}\big(F(X_{j-1},X_{j})X_{j}^{\prime}\theta-b(X_{j}^{\prime}\theta)\big) \\
&\qquad\qquad-\int\big(F(x_{1},x_{2})x_{2}^{\prime}\theta-b(x_{2}^{\prime}\theta)\big)P(x_{1},dx_{2})\nu(dx_{1})\cip0,
\end{align*}
where $P$ is the transition function.
}\qed
\end{rem}

\begin{thm}
Under Assumptions \ref{Ass1}--\ref{Ass3}, the QMLE satisfies 
\begin{align*}
\tes\cip\tz
\end{align*}
as $n\to\infty$.
\label{consis}
\end{thm}

\begin{ass}
(i) $\{X_{j};j=1,2,\ldots\}$ is strictly stationary. \\

(ii) For some $\Sigma_{0}>0$,

$\hspace{6mm}\displaystyle \frac{1}{n}E\bigg[\bigg\{\sumj\big(Y_{j}-\p b(X_{j}^{\prime}\tz)\big)X_{j}\bigg\}\bigg\{\sumj\big(Y_{j}-\p b(X_{j}^{\prime}\tz)\big)X_{j}\bigg\}^{\prime}\bigg]\to\Sig_{0}$.
\label{Ass6}
\end{ass}

\begin{thm}
Under Assumptions \ref{Ass1}--\ref{Ass3} and \ref{Ass6}, the asymptotic distribution of the QMLE is normal:
\begin{align*}
\sqrt{n}(\tes-\tz)\cil N(0,\Gam_{0}^{-1}\Sig_{0}\Gam_{0}^{-1}).
\end{align*}
\label{AN}
\end{thm}

Theorems \ref{consis} and \ref{AN} are shown in Supplementary Material.

Assume that Assumptions \ref{Ass1}--\ref{Ass3} and \ref{Ass6} (i) are satisfied.
When the candidate model is correctly specified, $\Sigma_{0}=\Gam_{0}$, i.e., $\sqrt{n}(\tes-\tz)\cil N(0,\Gam_{0}^{-1})$.
This is because the correctly specified model gives the equations $E_{\tz}[Y_{j}|X_{j}]=\p b_{m}(X_{j}^{\prime}\tz)$, $V_{\tz}[Y_{j}|X_{j}]=\p^{2} b_{m}(X_{j}^{\prime}\tz)$, and
\begin{align*}
&\frac{1}{n}E\bigg[\bigg\{\sumj\big(Y_{j}-\p b(X_{j}^{\prime}\tz)\big)X_{j}\bigg\}\bigg\{\sumj\big(Y_{j}-\p b(X_{j}^{\prime}\tz)\big)X_{j}\bigg\}^{\prime}\bigg] \\
&=\frac{1}{n}\sumj E\Big[E\big[\big(Y_{j}-\p b(X_{j}^{\prime}\tz)\big)^{2}|\mcf_{j-1}\vee\sigma(X_{j})\big]X_{j}X_{j}^{\prime}\Big] \\
&\qquad+2\sum_{i<j}\frac{1}{n}E\Big[\big(Y_{i}-\p b(X_{i}^{\prime}\tz)\big)E\big[\big(Y_{j}-\p b(X_{j}^{\prime}\tz)\big)|\mcf_{j-1}\vee\sigma(X_{j})\big]X_{i}X_{j}^{\prime}\Big] \\
&=\frac{1}{n}\sumj E\Big[E_{\tz}\big[\big(Y_{j}-E_{\tz}[Y_{j}|X_{j}]\big)^{2}|X_{j}\big]X_{j}X_{j}^{\prime}\Big] \\
&\qquad+2\sum_{i<j}\frac{1}{n}E\Big[\big(Y_{i}-\p b(X_{i}^{\prime}\tz)\big)E_{\tz}\big[\big(Y_{j}-E_{\tz}[Y_{j}|X_{j}]\big)|X_{j}\big]X_{i}X_{j}^{\prime}\Big] \\
&=\frac{1}{n}\sumj E\Big[\p^{2}b(X_{j}^{\prime}\tz)X_{j}X_{j}^{\prime}\Big]+0 \\
&=\int\p^{2}b(x^{\prime}\tz)xx^{\prime}\nu(dx).
\end{align*}

\begin{rem}{\rm
Suppose that Assumptions \ref{Ass2}--\ref{Ass3} and \ref{Ass6} (i) hold. 
The condition 
\begin{align*}
\left|\frac{1}{n}\sumj E\big[\big(Y_{j}-\p b(X_{j}^{\prime}\tz)\big)^{2}X_{j}X_{j}^{\prime}\big]-\Sigma_{0}\right|\to0
\end{align*}
implies Assumption \ref{Ass6} (ii), since it follows from this condition and Doukhan \cite[Theorem 3]{Dou94} that we obtain
\begin{align*}
&\left|\frac{1}{n}E\bigg[\bigg\{\sumj\big(Y_{j}-\p b(X_{j}^{\prime}\tz)\big)X_{j}\bigg\}\bigg\{\sumj\big(Y_{j}-\p b(X_{j}^{\prime}\tz)\big)X_{j}\bigg\}^{\prime}\bigg]-\Sigma_{0}\right| \\
&\leq\left|\frac{1}{n}\sumj E\big[\big(Y_{j}-\p b(X_{j}^{\prime}\tz)\big)^{2}X_{j}X_{j}^{\prime}\big]-\Sigma_{0}\right| \\
&\qquad+\frac{2}{n}\sum_{i<j}\Big|\cov\big[\big(Y_{i}-\p b(X_{i}^{\prime}\tz)\big)X_{i},(Y_{j}-\p b(X_{j}^{\prime}\tz)\big)X_{j}\big]\Big| \\
&\lesssim\left|\frac{1}{n}\sumj E\big[\big(Y_{j}-\p b(X_{j}^{\prime}\tz)\big)^{2}X_{j}X_{j}^{\prime}\big]-\Sigma_{0}\right|+\frac{2}{n}\sum_{i<j}\alpha(j-i)\to0
\end{align*}
as $n\to\infty$.
}\qed
\end{rem}


\section{Quasi-Bayesian information criterion for dependent GLM} \label{qbic.sec}
\subsection{Stochastic expansion}


\begin{ass}
$\displaystyle\frac{1}{\sqrt{n}}\sumj\big(F(X_{j})-\p b(X_{j}^{\prime}\tz)\big)X_{j}=O_{p}(1)$.
\label{Ass7}
\end{ass}

The next theorem shows the asymptotic behavior of the log marginal quasi-likelihood function.

\begin{thm}
Assume that Assumptions \ref{Ass1}--\ref{Ass3} and \ref{Ass7} hold and that the following conditions are satisfied:
\begin{itemize}
\item[(i)] $\displaystyle\pi(\tz)>0$, $\displaystyle\sup_{\theta\in\Theta}\pi(\theta)<\infty$.
\item[(ii)] For every $M>0$, $\displaystyle\sup_{|u|<M}\bigg|\pi\bigg(\tz+\frac{u}{\sqrt{n}}\bigg)-\pi(\tz)\bigg|\to0$ as $n\to\infty$.
\item[(iii)] $\displaystyle\log\pi(\tes)-\log\pi(\tz)=o_{p}(1)$.
\end{itemize}
Then, we have the expansion
\begin{align*}
&\log\bigg(\int_{\Theta}\exp\{\mbbh_{n}(\theta)\}\pi(\theta)d\theta\bigg) \\
&=\sumj\big(Y_{j}X_{j}^{\prime}\tes-b(X_{j}^{\prime}\tes)\big)-\frac{p}{2}\log n+\frac{p}{2}\log2\pi \\
&\qquad\qquad-\frac{1}{2}\log\det\bigg(\frac{1}{n}\sumj\p^{2}b(X_{j}^{\prime}\tes)X_{j}X_{j}^{\prime}\bigg)+\log\pi(\tes)+o_{p}(1) \\
&=\sumj\big(Y_{j}X_{j}^{\prime}\tes-b(X_{j}^{\prime}\tes)\big)+\frac{p}{2}\log2\pi \\
&\qquad\qquad-\frac{1}{2}\log\det\bigg(\sumj\p^{2}b(X_{j}^{\prime}\tes)X_{j}X_{j}^{\prime}\bigg)+\log\pi(\tes)+o_{p}(1).
\end{align*}
\label{qlfTh1}
\end{thm}

\begin{rem}{\rm
Suppose that we replace Assumptions \ref{Ass4} and \ref{Ass3} (ii) with the following conditions: 
\smallskip

(i) $\displaystyle\frac{1}{n}\sumj\big(F(X_{j})X_{j}^{\prime}\theta-b(X_{j}^{\prime}\theta)\big)\to\int\big(F(x)x^{\prime}\theta-b(x^{\prime}\theta)\big)\nu(dx)$ almost surely as $n\to\infty$, uniformly in $\theta\in\Theta$.
\smallskip

(ii) $\displaystyle\frac{1}{n}\sumj\p^{2}b(X_{j}^{\prime}\theta)X_{j}X_{j}^{\prime}\to\int\p^{2}b(x^{\prime}\theta)xx^{\prime}\nu(dx)$ almost surely as $n\to\infty$, uniformly in $\theta\in\Theta$.
\smallskip

(iii) For some constant $\lambda_{0}>0$, $\displaystyle P\bigg[\limsup_{n\to\infty}\lambda_{\min}\bigg(\frac{1}{n}\sumj\underline{b}(X_{j})X_{j}X_{j}^{\prime}\bigg)<\lambda_{0}\bigg]=0$. \\
Then, we can show that the log marginal quasi-likelihood function almost surely satisfies the expansion similar to Theorem \ref{qlfTh1}, i.e., almost surely
\begin{align*}
\log\bigg(\int_{\Theta}\exp\{\mbbh_{n}(\theta)\}\pi(\theta)d\theta\bigg)&=\sumj\big(Y_{j}X_{j}^{\prime}\tes-b(X_{j}^{\prime}\tes)\big)-\frac{p}{2}\log n+\frac{p}{2}\log2\pi \\
&\qquad\quad-\frac{1}{2}\log\det\bigg(\frac{1}{n}\sumj\p^{2}b(X_{j}^{\prime}\tes)X_{j}X_{j}^{\prime}\bigg) \\
&\qquad\quad+\log\pi(\tes)+o(1).
\end{align*}
}\qed
\end{rem}

Due to Theorem \ref{qlfTh1}, we define the quasi-Bayesian information criterion (QBIC) and BIC by
\begin{align*}
\qbic&=-2\sumj\big(Y_{j}X_{j}^{\prime}\tes-b(X_{j}^{\prime}\tes)\big)+\log\det\bigg(\sumj\p^{2}b(X_{j}^{\prime}\tes)X_{j}X_{j}^{\prime}\bigg)
\end{align*}
and
\begin{align*}
\mathrm{BIC}&=-2\sumj\big(Y_{j}X_{j}^{\prime}\tes-b(X_{j}^{\prime}\tes)\big)+p\log n,
\end{align*}
respectively.
Let $\qbic^{(1)},\ldots,\qbic^{(M)}$ be the QBIC for each candidate model. 
We calculate $\qbic^{(1)},\ldots,$ $\qbic^{(M)}$ and select the best model $\mfM_{m_{0}}$ having the minimum-QBIC value:
\begin{align*}
m_{0}=\argmin_{m\in\{1,\ldots,M\}}\qbic^{(m)}.
\end{align*}
The best model can be selected using BIC in a similar manner.
As directly seen by the definition, the QBIC have more computational load than the BIC.
Since the QBIC involves the observed-information matrix quantity, which is directly computed from data, the QBIC would more effectively take data dependence into account.
Furthermore, the penalty (second-term) of the QBIC consists of the second derivative of $\mbbh_{n}$ with respect to $\theta$, so the dimension of the parameter affects the penalty. 
That is, the QBIC implicitly takes the complexity of the model into account.


\subsection{Model selection consistency}

This section is devoted to verifying the model selection consistency of the (Q)BIC.
Let $\Theta_{i}\subset\mbbr^{p_{i}}$ and $\Theta_{j}\subset\mbbr^{p_{j}}$ be the parameter space associated with $\mfM_{i}$ and $\mfM_{j}$, respectively. 
If  $p_{i}<p_{j}$ and there exist a matrix $A\in\mbbr^{p_{j}\times p_{i}}$ with $A^{\prime}A=I_{p_{i}\times p_{i}}$ as well as a $c\in\mbbr^{p_{j}}$ such that $\mbbh_{i,n}(\theta)=\mbbh_{j,n}(A\theta+c)$ for all $\theta\in\Theta_{i}$, we say that $\Theta_{i}$ is nested in $\Theta_{j}$.
That is, when $\Theta_{i}$ is nested in $\Theta_{j}$, $\mfM_{j}$ includes $\mfM_{i}$.

Under Assumptions \ref{Ass1}--\ref{Ass3}, when $m_{0}$ satisfies 
\begin{align*}
\{m_{0}\}=\argmin_{m\in\mathcal{M}}\dim(\Theta_{m}),
\end{align*}
where $\mathcal{M}=\argmax_{m\in\{1,\ldots,M\}}\mbbh_{m,0}(\theta_{m,0})=\argmax_{m\in\{1,\ldots,M\}}\int\big(F(x)x^{\prime}\theta_{m,0}-b_{m}(x^{\prime}\theta_{m,0})\big)\nu(dx)$, we say that $\mfM_{m_{0}}$ is the {\it optimal} model.

\begin{thm}
Assume that Assumptions \ref{Ass1}--\ref{Ass3} and \ref{Ass7} are satisfied and that there exists a unique $m_{0}\in\{1,\ldots,M\}$ such that $\mfM_{m_{0}}$ is the optimal model. 
For any fixed $m\in\{1,\ldots,M\}\backslash\{m_{0}\}$, if $\Theta_{m_{0}}$ is nested in $\Theta_{m}$, or $\mbbh_{m,0}(\theta)\neq\mbbh_{m_{0},0}(\theta_{m_{0},0})$ for any $\theta\in\Theta_{m}$, then 
\begin{align*}
\lim_{n\to\infty}P\big[\qbic^{(m_{0})}-\qbic^{(m)}<0\big]=1.
\end{align*}
\label{mod_consis}
\end{thm}

This theorem implies that the probability that the optimal model is selected by using QBIC tends to 1 as $n\to\infty$.
The probability that BIC choose the optimal model can be handled analogously.


\section{Examples and simulation results} \label{ExSim}
This section presents the results of simulations to evaluate the finite sample performance of the model selection using QBIC, BIC, and formal AIC (fAIC). 
Because the theoretical part of AIC is not dealt with, we use the word fAIC as the AIC, i.e., the fAIC of the $m$th model is defined by
\begin{align*}
\mathrm{fAIC}^{(m)}=-2\mbbh_{m,n}(\hat{\theta}_{m,n})+2p_{m}.
\end{align*}



\subsection{Model selection in a correctly specified model} \label{Simu1}
We assume that the explanatory variables $X_{j,1},...,X_{j,4}$ are given by 
\begin{align*}
& X_{j,1}=1\;(j\geq1), \\
& X_{1,2}=1,\;X_{j,2}=0.5X_{j-1,2}+\epsilon_{j,2},\;(j\geq2), \\
& X_{1,3}=0,\;X_{j,3}=-0.7X_{j-1,3}+\epsilon_{j,3},\;(j\geq2), \\
& X_{1,4}=-1,\;X_{j,4}=0.8X_{j-1,4}+\epsilon_{j,4},\;(j\geq2),
\end{align*}
where the error vector $(\epsilon_{j,2},\epsilon_{j,3},\epsilon_{j,4})\sim N(0,\Sig)$ with $\Sig=(0.5^{|k-\ell|})_{k,\ell=1,2,3}$.
Moreover, the response variable $Y_{j}$ is obtained from the true model defined by the linear logistic regression model
\begin{align}
Y_{j}\sim B\bigg(1,\frac{\exp(X_{j}^{\prime}\theta^{\ast})}{1+\exp(X_{j}^{\prime}\theta^{\ast})}\bigg),
\label{simu1t}
\end{align}
where the true value $\theta^{\ast}=(0,-3,0,1)$ and $B(1, P)$ is a Bernoulli distribution with success probability $P$. 
The candidate models are given by the model (\ref{simu1t}) and consist of the following combination of $X_{j}$:
\begin{align*}
&{\bf Model}\;{\bf1:}\;X_{j}=(X_{j,1},X_{j,2},X_{j,3},X_{j,4});\;{\bf Model}\;{\bf2:}\;X_{j}=(X_{j,1},X_{j,2},X_{j,3}); \\
&{\bf Model}\;{\bf3:}\;X_{j}=(X_{j,1},X_{j,2},X_{j,4});\;{\bf Model}\;{\bf4:}\;X_{j}=(X_{j,1},X_{j,3},X_{j,4}); \\
&{\bf Model}\;{\bf5:}\;X_{j}=(X_{j,2},X_{j,3},X_{j,4});\;{\bf Model}\;{\bf6:}\;X_{j}=(X_{j,1},X_{j,2}); \\
&{\bf Model}\;{\bf7:}\;X_{j}=(X_{j,1},X_{j,3});\;{\bf Model}\;{\bf8:}\;X_{j}=(X_{j,1},X_{j,4}); \\
&{\bf Model}\;{\bf9:}\;X_{j}=(X_{j,2},X_{j,3});\;\;{\bf Model}\;{\bf10:}\;X_{j}=(X_{j,2},X_{j,4}); \\
&{\bf Model}\;{\bf11:}\;X_{j}=(X_{j,3},X_{j,4});\;{\bf Model}\;{\bf12:}\;X_{j}=X_{j,1};\;{\bf Model}\;{\bf13:}\;X_{j}=X_{j,2}; \\
&{\bf Model}\;{\bf14:}\;X_{j}=X_{j,3};\;{\bf Model}\;{\bf15:}\;X_{j}=X_{j,4}.
\end{align*}
Then, the optimal model is Model 10, and Models 1, 3, 5 contain the optimal model. 
The number of models selected using QBIC, BIC, and fAIC from among Models 1--15 over 10,000 simulations was counted. 
For example, in the case of Model 1, $\mbbh_{1,n}$, QBIC, BIC, and fAIC are given by
\begin{align*}
\mbbh_{1,n}(\theta)&=\sumj\bigg\{Y_{j}\sum_{i=1}^{4}X_{j,i}\theta_{i}-\log\bigg(1+\exp\Big(\sum_{i=1}^{4}X_{j,i}\theta_{i}\Big)\bigg)\bigg\}, \\
\qbic&=-2\mbbh_{1,n}(\tes)+\log\det\left(\sumj\frac{\exp\Big(\sum_{i=1}^{4}X_{j,i}\hat{\theta}_{i,n}\Big)X_{j}X_{j}^{\prime}}{\Big(1+\exp\Big(\sum_{i=1}^{4}X_{j,i}\hat{\theta}_{i,n}\Big)\Big)^{2}}\right), \\
\mathrm{BIC}&=-2\mbbh_{1,n}(\tes)+4\log n
\end{align*}
and
\begin{align*}
\mathrm{fAIC}&=-2\mbbh_{1,n}(\tes)+4\times2,
\end{align*}
where the QMLE $\tes=(\hat{\theta}_{1,n},\hat{\theta}_{2,n},\hat{\theta}_{3,n},\hat{\theta}_{4,n})$ maximizes $\mbbh_{1,n}$.
For numerical optimization, we set the initial values to be random numbers generated from uniform distribution $U(\theta^{\ast}-1,\,\theta^{\ast}+1)$.

Let us verify the assumptions for (Q)BIC.
In the current case, the function $b$ defined in (\ref{glm_pdf}) is given by $b(\theta)=\log(1+e^{\theta})$, 
and Assumptions \ref{Ass1} and \ref{Ass2} are satisfied with $E[Y_{j}|\mcf_{j-1}\vee\sigma(X_{j})]=F(X_{j})=e^{-3X_{j,2}+X_{j,4}}/(1+e^{-3X_{j,2}+X_{j,4}})$. 
In particular, $ \sup_{j\in\mbbn}E\big[|X_{j}|^{q}\big]<\infty$ for every $q>0$.
Furthermore, the function $\underline{b}$ of Assumption \ref{Ass3} can be given by $\underline{b}(x)=e^{-C|x|}/(1+e^{C|x|})^{2}$ for some constant $C>0$ satisfying $\sup_{\theta\in\Theta}|x^{\prime}\theta|\leq C|x|$.
If $\{(X_{j},Y_{j});j=1,2,\ldots\}$ has the exponential $\beta$-mixing property, then Assumptions \ref{Ass5}, \ref{Ass4}, and \ref{Ass7} hold.
The sufficient conditions for the exponential $\beta$-mixing property were given by Baraud {\it et al.} \cite{BarComVie01} and Doukhan \cite[Section 2.4]{Dou94}.

Table \ref{select1} summarizes the comparison results of the model selection frequency.
Model 10 is selected with high frequency for all criteria and $n$. 
Moreover, the probability that Model 10 is selected by (Q)BIC increases as $n$ increases.
In Table \ref{esti1}, the differences between the true values and the estimators in the specified models decrease when $n$ increases.
These results demonstrate the consistency of the estimators and the model selection consistency of QBIC and BIC.

\begin{rem}{\rm
In misspecified models, it may happen that optimal parameter $\tz\ne\theta^{\ast}$.
Then, the estimators are not necessarily estimating the true values.
}\qed
\end{rem}

\begin{table}
\caption{\footnotesize The numbers of models selected by QBIC, BIC, and fAIC in Section \ref{Simu1} over 10,000 simulations for various $n$ (1--15 represent the models, and the optimal model is Model 10)}
\scalebox{0.83}[0.8]{
\begin{tabular}{ l c c c c c c c c c c c c c c c } \hline
\multicolumn{1}{l}{Criteria} & \multicolumn{15}{c}{$n=50$} \\
 & 1 & 2 & 3 & 4 & 5 & 6 & 7 & 8 & 9 & $10^{\ast}$ & 11 & 12 & 13 & 14 & 15 \\ \hline
QBIC & 1489 & 65 & 2084 & 0 & 1260 & 201 & 0 & 0 & 56 & {\bf 4666} & 0 & 0 & 172 & 0 & 0 \\ 
BIC & 99 & 18 & 531 & 0 & 562 & 222 & 0 & 0 & 86 & {\bf 7720} & 0 & 0 & 762 & 0 & 0 \\ 
fAIC & 479 & 55 & 1310 & 0 & 1424 & 192 & 0 & 0 & 93 & {\bf 6242} & 0 & 0 & 205 & 0 & 0 \\ \hline
\multicolumn{1}{l}{Criteria} & \multicolumn{15}{c}{$n=100$} \\
 & 1 & 2 & 3 & 4 & 5 & 6 & 7 & 8 & 9 & $10^{\ast}$ & 11 & 12 & 13 & 14 & 15 \\ \hline
QBIC & 298 & 2 & 1483 & 0 & 989 & 10 & 0 & 0 & 4 & {\bf 7206} & 0 & 0 & 8 & 0 & 0 \\ 
BIC & 19 & 1 & 323 & 0 & 397 & 15 & 0 & 0 & 8 & {\bf 9179} & 0 & 0 & 58 & 0 & 0 \\ 
fAIC & 347 & 1 & 1380 & 0 & 1367 & 7 & 0 & 0 & 2 & {\bf 6895} & 0 & 0 & 1 & 0 & 0 \\ \hline
\multicolumn{1}{l}{Criteria} & \multicolumn{15}{c}{$n=200$} \\
 & 1 & 2 & 3 & 4 & 5 & 6 & 7 & 8 & 9 & $10^{\ast}$ & 11 & 12 & 13 & 14 & 15 \\ \hline
QBIC & 86 & 0 & 910 & 0 & 616 & 0 & 0 & 0 & 0 & {\bf 8388} & 0 & 0 & 0 & 0 & 0 \\ 
BIC & 5 & 0 & 235 & 0 & 222 & 0 & 0 & 0 & 0 & {\bf 9538} & 0 & 0 & 0 & 0 & 0 \\ 
fAIC & 281 & 0 & 1314 & 0 & 1414 & 0 & 0 & 0 & 0 & {\bf 6991} & 0 & 0 & 0 & 0 & 0 \\ \hline 
\end{tabular}
}
\label{select1}
\end{table}

\begin{table}
\caption{\footnotesize The mean and standard deviation (s.d.) of the estimators $\hat{\theta}_{1}$, $\hat{\theta}_{2}$, $\hat{\theta}_{3}$, and $\hat{\theta}_{4}$ for various $n$ (1--15 represent the models, and the true value $\theta^{\ast}=(0,-3,0,1)$)}
\scalebox{0.56}[0.7]{
\begin{tabular}{ l l c c c c c c c c c c c c c c } \hline
\multicolumn{1}{l}{} & \multicolumn{1}{l}{}  & \multicolumn{4}{c}{$n=50$} & \multicolumn{1}{l}{} & \multicolumn{4}{c}{$n=100$} & \multicolumn{1}{l}{} & \multicolumn{4}{c}{$n=200$} \\ \cline{3-6} \cline{8-11} \cline{13-16}
\multicolumn{16}{l}{} \\[-10pt]
 & & $\hat{\theta}_{1}$ & $\hat{\theta}_{2}$ & $\hat{\theta}_{3}$ & $\hat{\theta}_{4}$ & & $\hat{\theta}_{1}$ & $\hat{\theta}_{2}$ & $\hat{\theta}_{3}$ & $\hat{\theta}_{4}$ & & $\hat{\theta}_{1}$ & $\hat{\theta}_{2}$ & $\hat{\theta}_{3}$ & $\hat{\theta}_{4}$ \\ \hline
1 & mean & -0.0793 & -8.3409 & -0.0219 & 2.7990 & & 0.0004 & -3.3727 & -0.0057 & 1.1266 & & 0.0023 & -3.1642 & 0.0004 & 1.0542 \\ 
 & s.d. & 8.1943 & 35.1378 & 5.9571 & 12.8200 & & 0.3895 & 0.8918 & 0.2889 & 0.3744 & & 0.2425 & 0.4946 & 0.1807 & 0.2123 \\ \hline 
2 & mean & -0.0505 & -2.7481 & 0.1061 & -- & & -0.0425 & -2.1878 & 0.0867 & -- & & -0.0167 & -2.0734 & 0.0856 & -- \\ 
 & s.d. & 2.0544 & 7.9339 & 1.3571 & -- & & 0.4653 & 0.5332 & 0.2141 & -- & & 0.3146 & 0.3336 & 0.1395 & -- \\ \hline
3 & mean & 0.0355 & -5.7372 & -- & 1.8913 & & 0.0001 & -3.2941 & -- & 1.0993 & & 0.0021 & -3.1332 & -- & 1.0441 \\ 
 & s.d. & 4.6176 & 23.3197 & -- & 8.1816 & & 0.3763 & 0.8210 & -- & 0.3508 & & 0.2397 & 0.4804 & -- & 0.2078 \\ \hline
4 & mean & -0.0999 & -- & -0.2581 & 0.3139 & & -0.0451 & -- & -0.2364 & 0.2852 & & -0.0168 & -- & -0.2250 & 0.2746 \\ 
 & s.d. & 0.4791 & -- & 0.2318 & 0.3193 & & 0.3109 & -- & 0.1453 & 0.1940 & & 0.2100 & -- & 0.0976 & 0.1264 \\ \hline
5 & mean & -- & -5.6381 & 0.0219 & 1.8650 & & -- & -3.2928 & -0.0052 & 1.0989 & & -- & -3.1336 & 0.0003 & 1.0442 \\ 
 & s.d. & -- & 22.5995 & 3.7461 & 7.8695 & & -- & 0.8333 & 0.2792 & 0.3465 & & -- & 0.4835 & 0.1787 & 0.2057 \\ \hline
6 & mean & -0.0635 & -2.4621 & -- & -- & & -0.0429 & -2.1293 & -- & -- & & -0.0169 & -2.0324 & -- & -- \\ 
 & s.d. & 1.8364 & 4.1725 & -- & -- & & 0.4578 & 0.5074 & -- & -- & & 0.3127 & 0.3232 & -- & -- \\ \hline
7 & mean & -0.1086 & -- & -0.1960 & -- & & -0.0518 & -- & -0.1808 & -- & & -0.0217 & -- & -0.1723 & -- \\ 
 & s.d. & 0.4591 & -- & 0.2093 & -- & & 0.3188 & -- & 0.1351 & -- & & 0.2227 & -- & 0.0921 & -- \\ \hline
8 & mean & -0.1006 & -- & -- & 0.2681 & & -0.0453 & -- & -- & 0.2483 & & -0.0170 & -- & -- & 0.2415 \\ 
 & s.d. & 0.4660 & -- & -- & 0.3022 & & 0.3063 & -- & -- & 0.1875 & & 0.2081 & -- & -- & 0.1230 \\ \hline
9 & mean & -- & -2.3773 & 0.1058 & -- & & -- & -2.1041 & 0.0878 & -- & & -- & -2.0342 & 0.0860 & -- \\ 
 & s.d. & -- & 5.2112 & 1.0538 & -- & & -- & 0.5008 & 0.2061 & -- & & -- & 0.3263 & 0.1372 & -- \\ \hline
$10^{\ast}$ & mean & -- & -4.2068 & -- & 1.3952 & & -- & -3.2211 & -- & 1.0741 & & -- & -3.1037 & -- & 1.0344 \\ 
 & s.d. & -- & 13.1535 & -- & 4.3787 & & -- & 0.7702 & -- & 0.3259 & & -- & 0.4699 & -- & 0.2013 \\ \hline
11 & mean & -- & -- & -0.2546 & 0.3124 & & -- & -- & -0.2350 & 0.2855 & & -- & -- & -0.2243 & 0.2747 \\ 
 & s.d. & -- & -- & 0.2218 & 0.2851 & & -- & -- & 0.1424 & 0.1840 & & -- & -- & 0.0967 & 0.1230 \\ \hline
12 & mean & -0.0695 & -- & -- & -- & & -0.0368 & -- & -- & -- & & -0.0179 & -- & -- & -- \\ 
 & s.d. & 0.3186 & -- & -- & -- & & 0.2548 & -- & -- & -- & & 0.1936 & -- & -- & -- \\ \hline
13 & mean & -- & -2.6475 & -- & -- & & -- & -2.5688 & -- & -- & & -- & -2.5278 & -- & -- \\ 
 & s.d. & -- & 0.3694 & -- & -- & & -- & 0.3052 & -- & -- & & -- & 0.2789 & -- & -- \\ \hline
14 & mean & -- & -- & -0.1525 & -- & & -- & -- & -0.1528 & -- & & -- & -- & -0.1517 & -- \\ 
 & s.d. & -- & -- & 0.1694 & -- & -- & & -- & 0.1218 & -- & & -- & -- & 0.0901 & -- \\ \hline
15 & mean & -- & -- & -- & 0.5703 & & -- & -- & -- & 0.5507 & & -- & -- & -- & 0.5394 \\ 
 & s.d. & -- & -- & -- & 0.2553 & & -- & -- & -- & 0.2478 & & -- & -- & -- & 0.2477 \\ \hline
\end{tabular}
}
\label{esti1}
\end{table}


\subsection{Model selection in a misspecified model} \label{Simu2}
We use the same conditions as in the previous section except that the response variable $Y_{j}$ is obtained from the true model defined by 
\begin{align*}
Y_{j}\sim B\Big(1,\Phi(X_{j}^{\prime}\theta^{\ast})\Big),
\end{align*}
where $\Phi(x)=\int_{-\infty}^{x}\frac{1}{\sqrt{2\pi}}\exp(-\frac{t^{2}}{2})dt$. 
Then, Models 1--15 are misspecified models.
As with Section \ref{Simu1}, we can show the validity of the assumptions for (Q)BIC.

From Table \ref{select2}, we obtain similar results even though the candidate models do not include the true model.
Table \ref{esti2} summarizes the mean and standard deviation of the estimators. 
Since the optimal parameter values are not given here, we can not see the differences between the optimal parameter values and the estimators, although the standard deviations decrease as $n$ increases.

\begin{table}
\caption{\footnotesize The numbers of models selected by QBIC, BIC, and fAIC in Section \ref{Simu2} over 10,000 simulations for various $n$ (1--15 represent the models)}
\scalebox{0.9}[0.8]{
\begin{tabular}{ l c c c c c c c c c c c c c c c } \hline
\multicolumn{1}{l}{Criteria} & \multicolumn{15}{c}{$n=100$} \\
 & 1 & 2 & 3 & 4 & 5 & 6 & 7 & 8 & 9 & 10 & 11 & 12 & 13 & 14 & 15 \\ \hline
QBIC & 965 & 0 & 2025 & 0 & 1321 & 0 & 0 & 0 & 0 & 5689 & 0 & 0 & 0 & 0 & 0 \\ 
BIC & 41 & 0 & 435 & 0 & 398 & 0 & 0 & 0 & 0 & 9125 & 0 & 0 & 1 & 0 & 0 \\ 
fAIC & 443 & 0 & 1452 & 0 & 1538 & 0 & 0 & 0 & 0 & 6567 & 0 & 0 & 0 & 0 & 0 \\  \hline
\multicolumn{1}{l}{Criteria} & \multicolumn{15}{c}{$n=200$} \\
 & 1 & 2 & 3 & 4 & 5 & 6 & 7 & 8 & 9 & 10 & 11 & 12 & 13 & 14 & 15 \\ \hline
QBIC & 223 & 0 & 1338 & 0 & 915 & 0 & 0 & 0 & 0 & 7524 & 0 & 0 & 0 & 0 & 0 \\ 
BIC & 9 & 0 & 278 & 0 & 274 & 0 & 0 & 0 & 0 & 9439 & 0 & 0 & 0 & 0 & 0 \\ 
fAIC & 349 & 0 & 1436 & 0 & 1414 & 0 & 0 & 0 & 0 & 6801 & 0 & 0 & 0 & 0 & 0 \\ \hline
\multicolumn{1}{l}{Criteria} & \multicolumn{15}{c}{$n=300$} \\
 & 1 & 2 & 3 & 4 & 5 & 6 & 7 & 8 & 9 & 10 & 11 & 12 & 13 & 14 & 15 \\ \hline
QBIC & 108 & 0 & 1009 & 0 & 694 & 0 & 0 & 0 & 0 & 8189 & 0 & 0 & 0 & 0 & 0 \\ 
BIC & 5 & 0 & 190 & 0 & 216 & 0 & 0 & 0 & 0 & 9589 & 0 & 0 & 0 & 0 & 0 \\ 
fAIC & 295 & 0 & 1352 & 0 & 1388 & 0 & 0 & 0 & 0 & 6965 & 0 & 0 & 0 & 0 & 0 \\ \hline 
\end{tabular}
}
\label{select2}
\end{table}

\begin{table}
\caption{\footnotesize The mean and standard deviation (s.d.) of the estimators $\hat{\theta}_{1}$, $\hat{\theta}_{2}$, $\hat{\theta}_{3}$, and $\hat{\theta}_{4}$ in each model for various $n$ (1--15 represent the models)}
\scalebox{0.58}[0.7]{
\begin{tabular}{ l l c c c c c c c c c c c c c c } \hline
\multicolumn{1}{l}{} & \multicolumn{1}{l}{}  & \multicolumn{4}{c}{$n=100$} & \multicolumn{1}{l}{} & \multicolumn{4}{c}{$n=300$} & \multicolumn{1}{l}{} & \multicolumn{4}{c}{$n=300$} \\ \cline{3-6} \cline{8-11} \cline{13-16}
\multicolumn{16}{l}{} \\[-10pt]
 & & $\hat{\theta}_{1}$ & $\hat{\theta}_{2}$ & $\hat{\theta}_{3}$ & $\hat{\theta}_{4}$ & & $\hat{\theta}_{1}$ & $\hat{\theta}_{2}$ & $\hat{\theta}_{3}$ & $\hat{\theta}_{4}$ & & $\hat{\theta}_{1}$ & $\hat{\theta}_{2}$ & $\hat{\theta}_{3}$ & $\hat{\theta}_{4}$ \\ \hline
1 & mean & 0.0020 & 7.6472 & -0.0028 & -2.5525 & & -0.0031 & 5.8988 & -0.0009 & -1.9640 & & 0.0038 & 5.6861 & -0.0006 & -1.8975 \\ 
 & s.d. & 1.7319 & 14.4366 & 1.2513 & 5.0473 & & 0.3440 & 1.1541 & 0.2541 & 0.4378 & & 0.2599 & 0.8285 & 0.1914 & 0.3118 \\ \hline
2 & mean & 0.0596 & 2.8857 & -0.1171 & -- & & 0.0255 & 2.6805 & -0.1079 & -- & & 0.0176 & 2.6080 & -0.1066 & -- \\ 
 & s.d. & 0.5917 & 0.7518 & 0.2335 & -- & & 0.3873 & 0.4610 & 0.1501 & -- & & 0.3069 & 0.3574 & 0.1168 & -- \\ \hline
3 & mean & 0.0104 & 6.7840 & -- & -2.2605 & & -0.0030 & 5.7925 & -- & -1.9290 & & 0.0036 & 5.6262 & -- & -1.8772 \\ 
 & s.d. & 1.0092 & 9.8524 & -- & 3.4170 & & 0.3351 & 1.0856 & -- & 0.4150 & & 0.2561 & 0.8056 & -- & 0.3038 \\ \hline
4 & mean & 0.0525 & -- & 0.2555 & -0.3135 & & 0.0218 & -- & 0.2448 & -0.2959 & & 0.0157 & -- & 0.2415 & -0.2931 \\ 
 & s.d. & 0.3315 & -- & 0.1409 & 0.2024 & & 0.2209 & -- & 0.0938 & 0.1314 & & 0.1781 & -- & 0.0756 & 0.1048 \\ \hline
5 & mean & -- & 6.7809 & 0.0061 & -2.2594 & & -- & 5.7915 & -0.0008 & -1.9285 & & -- & 5.6255 & -0.0007 & -1.8771 \\ 
 & s.d. & -- & 9.5017 & 0.8186 & 3.3869 & & -- & 1.0916 & 0.2479 & 0.4130 & & -- & 0.8008 & 0.1886 & 0.3009 \\ \hline
6 & mean & 0.0594 & 2.7967 & -- & -- & & 0.0256 & 2.6226 & -- & -- & & 0.0177 & 2.5586 & -- & -- \\ 
 & s.d. & 0.5801 & 0.7108 & -- & -- & & 0.3839 & 0.4437 & -- & -- & & 0.3053 & 0.3487 & -- & -- \\ \hline
7 & mean & 0.0619 & -- & 0.1933 & -- & & 0.0276 & -- & 0.1866 & -- & & 0.0189 & -- & 0.1838 & -- \\ 
 & s.d. & 0.3385 & -- & 0.1298 & -- & & 0.2324 & -- & 0.0877 & -- & & 0.1902 & -- & 0.0712 & -- \\ \hline
8 & mean & 0.0527 & -- & -- & -0.2728 & & 0.0220 & -- & -- & -0.2591 & & 0.0158 & -- & -- & -0.2574 \\ 
 & s.d. & 0.3260 & -- & -- & 0.1951 & & 0.2182 & -- & -- & 0.1276 & & 0.1761 & -- & -- & 0.1019 \\ \hline
9 & mean & -- & 2.7283 & -0.1163 & -- & & -- & 2.6104 & -0.1073 & -- & & -- & 2.5642 & -0.1064 & -- \\ 
 & s.d. & -- & 0.6923 & 0.2218 & -- & & -- & 0.4446 & 0.1466 & -- & & -- & 0.3506 & 0.1151 & -- \\ \hline
10 & mean & -- & 6.2763 & -- & -2.0919 & & -- & 5.6925 & -- & -1.8959 & & -- & 5.5676 & -- & -1.8575 \\ 
 & s.d. & -- & 6.7882 & -- & 2.4910 & & -- & 1.0309 & -- & 0.3934 & & -- & 0.7797 & -- & 0.2936 \\ \hline
11 & mean & -- & -- & 0.2535 & -0.3120 & & -- & -- & 0.2438 & -0.2954 & & -- & -- & 0.2409 & -0.2929 \\ 
 & s.d. & -- & -- & 0.1375 & 0.1922 & & -- & -- & 0.0927 & 0.1277 & & -- & -- & 0.0751 & 0.1028 \\ \hline
12 & mean & 0.0457 & -- & -- & -- & & 0.0218 & -- & -- & -- & & 0.0165 & -- & -- & -- \\ 
 & s.d. & 0.2639 & -- & -- & -- & & 0.1971 & -- & -- & -- & & 0.1677 & -- & -- & -- \\ \hline
13 & mean & -- & 2.7922 & -- & -- & & -- & 2.7190 & -- & -- & & -- & 2.6821 & -- & -- \\ 
 & s.d. & -- & 0.3799 & -- & -- & & -- & 0.3080 & -- & -- & & -- & 0.2692 & -- & -- \\ \hline
14 & mean & -- & -- & 0.1636 & -- & & -- & -- & 0.1630 & -- & & -- & -- & 0.1631 & -- \\ 
 & s.d. & -- & -- & 0.1182 & -- & & -- & -- & 0.0885 & -- & & -- & -- & 0.0756 & -- \\ \hline
15 & mean & -- & -- & -- & -0.5569 & & -- & -- & -- & -0.5398 & & -- & -- & -- & -0.5376 \\ 
 & s.d. & -- & -- & -- & 0.2432 & & -- & -- & -- & 0.2439 & & -- & -- & -- & 0.2466 \\ \hline
\end{tabular}
}
\label{esti2}
\end{table}


\subsection{Model selection in univariate time series model} \label{Simu3}
Let $X_{j}=(Z_{j},Z_{j-1},\ldots,Z_{j-(p-1)})^{\prime}$ be the explanatory vector for $j\in\{1,\ldots,n\}$, where for every $i\in\{2,\ldots,n\}$, $Z_{i}$ is given by
\begin{align*}
Z_{-n+2}=\cdots=Z_{0}=0,\;Z_{1}=1,\;Z_{i}=0.6Z_{i-1}+\epsilon_{i},
\end{align*}
where $\epsilon_{i}\sim N(0,1)$.
The response variable $Y_{j}$ is obtained from the true model defined by 
\begin{align}
Y_{j}\sim B\bigg(1,\frac{\exp(X_{j}^{\prime}\theta^{\ast})}{1+\exp(X_{j}^{\prime}\theta^{\ast})}\bigg),
\label{simu3t}
\end{align}
where the true value $\theta^{\ast}=(3,-1,2,1)$. 
For simplicity, we here focus on the hierarchical models as the candidate models:
\begin{align*}
&{\bf Model}\;{\bf1:}\;X_{j}=(Z_{j});\;{\bf Model}\;{\bf2:}\;X_{j}=(Z_{j},Z_{j-1}); \\
&{\bf Model}\;{\bf3:}\;X_{j}=(Z_{j},Z_{j-1},Z_{j-2});\;{\bf Model}\;{\bf4:}\;X_{j}=(Z_{j},Z_{j-1},Z_{j-2},Z_{j-3}); \\
&{\bf Model}\;{\bf5:}\;X_{j}=(Z_{j},Z_{j-1},Z_{j-2},Z_{j-3},Z_{j-4});\;\cdots.
\end{align*}
Then, the optimal model is Model 4. 

The number of models selected using QBIC, BIC, and fAIC from among the candidate models over 10,000 simulations were calculated. 
First, we calculate $\qbic^{(1)}$ and $\qbic^{(2)}$. 
If $\qbic^{(1)}<\qbic^{(2)}$, Model 1 is selected as the best model.
When $\qbic^{(1)}\geq\qbic^{(2)}$, we calculate $\qbic^{(3)}$ and compare $\qbic^{(2)}$ with $\qbic^{(3)}$.
The same procedures are repeated until they are stopped at the best model.
The cases of BIC and fAIC are calculated in a similar manner.
Note that the validity of the assumptions for (Q)BIC can be checked using a method similar to that described in Section \ref{Simu1}.

Table \ref{select3} summarizes the comparison results of the model selection frequency.
The best model is searched for among Models 1--12 for all cases.
Model 4 is selected as the best model with the highest frequency.
Moreover, the frequency that Model 4 is selected by (Q)BIC increases as $n$ increases.
This result demonstrates that QBIC and BIC have model selection consistency. 
In Table \ref{esti3}, the differences between the true values and the estimators in the specified models (Models 4--6) decrease as $n$ increases, and the standard deviations behave similarly.
Hence, the consistency of the estimators can be observed.

\begin{rem}{\rm
The current situation satisfies the original model setting given in Section \ref{NoTe} even if $\{Z_{j}; j=1,2,\ldots\}$ is replaced by the $m$th Markov chain ($m\geq 2$).
}\qed
\end{rem}

\begin{rem}{\rm
Since we only focus on the contribution of $\{Z_{j}\}$ to $\mbY_{n}$ in this simulation, the Bayesian model selection is possible without specifying the distribution of $\{Z_{j}\}$.
}\qed
\end{rem}

\begin{table}
\caption{\footnotesize The numbers of models selected by QBIC, BIC, and fAIC in Section \ref{Simu3} over 10,000 simulations for various $n$ (1--11 represent the models, and the optimal model is Model 4)}
\scalebox{0.93}[0.8]{
\begin{tabular}{ l c c c c c c c c c c c c } \hline
\multicolumn{1}{l}{Criteria} & \multicolumn{12}{c}{$n=100$} \\
 & 1 & 2 & 3 & $4^{\ast}$ & 5 & 6 & 7 & 8 & 9 & 10 & 11 & 12 \\ \hline
QBIC & 2811 & 0 & 621 & {\bf 4790} & 1230 & 389 & 116 & 39 & 2 & 1 & 1 & 0 \\ 
BIC & 4132 & 0 & 1787 & {\bf 3884} & 186 & 10 & 1 & 0 & 0 & 0 & 0 & 0 \\ 
fAIC & 1 & 0 & 539 & {\bf 6079} & 2220 & 832 & 235 & 84 & 18 & 1 & 0 & 1 \\ \hline
\multicolumn{1}{l}{Criteria} & \multicolumn{12}{c}{$n=200$} \\
 & 1 & 2 & 3 & $4^{\ast}$ & 5 & 6 & 7 & 8 & 9 & 10 & 11 & 12 \\ \hline
QBIC & 1412 & 0 & 137 & {\bf 7288} & 1012 & 130 & 19 & 1 & 1 & 0 & 0 & 0 \\ 
BIC & 2229 & 0 & 537 & {\bf 7059} & 170 & 3 & 2 & 0 & 0 & 0 & 0 & 0 \\ 
fAIC & 0 & 0 & 44 & {\bf 6601} & 2311 & 770 & 207 & 52 & 12 & 2 & 1 & 0 \\ \hline
\multicolumn{1}{l}{Criteria} & \multicolumn{12}{c}{$n=300$} \\
 & 1 & 2 & 3 & $4^{\ast}$ & 5 & 6 & 7 & 8 & 9 & 10 & 11 & 12 \\ \hline
QBIC & 788 & 0 & 24 & {\bf 8263} & 836 & 86 & 3 & 0 & 0 & 0 & 0 & 0 \\ 
BIC & 1252 & 0 & 131 & {\bf 8449} & 166 & 2 & 0 & 0 & 0 & 0 & 0 & 0 \\ 
fAIC & 0 & 0 & 3 & {\bf 6749} & 2329 & 695 & 167 & 49 & 6 & 2 & 0 & 0 \\ \hline 
\end{tabular}
}
\label{select3}
\end{table}

\begin{table}
\caption{\footnotesize The mean and standard deviation (s.d.) of the estimators $\hat{\theta}_{1}$, $\hat{\theta}_{2}$, $\hat{\theta}_{3}$, $\hat{\theta}_{4}$, $\hat{\theta}_{5}$, and $\hat{\theta}_{6}$ in each model for various $n$ (1--6 represent the models, and the true value $\theta^{\ast}=(3,-1,2,1)$)}
\scalebox{1.1}[0.67]{
\begin{tabular}{ l l  c c c c c c } \hline
\multicolumn{1}{l}{} & \multicolumn{1}{l}{} & \multicolumn{6}{c}{$n=100$}  \\
 & & $\hat{\theta}_{1}$ & $\hat{\theta}_{2}$ & $\hat{\theta}_{3}$ & $\hat{\theta}_{4}$ & $\hat{\theta}_{5}$ & $\hat{\theta}_{6}$ \\ \hline
1 & mean & 2.5172 & -- & -- & -- & -- & -- \\ 
 & s.d. & 0.2870 & -- & -- & -- & -- & -- \\ \hline 
2 & mean & 1.6640 & 0.2999 & -- & -- & -- & -- \\ 
 & s.d. & 0.4144 & 0.2984 & -- & -- & -- & -- \\ \hline 
3 & mean & 2.9344 & -0.9739 & 2.5388 & -- & -- & -- \\ 
 & s.d. & 1.6998 & 0.6616 & 1.6373 & -- & -- & -- \\ \hline 
$4^{\ast}$ & mean & 3.6868 & -1.2316 & 2.4792 & 1.2357 & -- & -- \\ 
 & s.d. & 4.7146 & 1.8259 & 3.7514 & 1.9396 & -- & -- \\ \hline 
5 & mean & 4.0140 & -1.3648 & 2.6861 & 1.3428 & 0.0084 & -- \\ 
 & s.d. & 7.1375 & 3.3559 & 5.1527 & 2.5983 & 1.3582 & -- \\ \hline 
6 & mean & 4.3051 & -1.4502 & 2.8596 & 1.4699 & -0.0248 & 0.0181 \\ 
 & s.d. & 9.2106 & 3.8772 & 6.0818 & 3.8938 & 2.1885 & 1.3463 \\ \hline 
\multicolumn{1}{l}{} & \multicolumn{1}{l}{} & \multicolumn{6}{c}{$n=200$} \\
  & & $\hat{\theta}_{1}$ & $\hat{\theta}_{2}$ & $\hat{\theta}_{3}$ & $\hat{\theta}_{4}$ & $\hat{\theta}_{5}$ & $\hat{\theta}_{6}$ \\ \hline
1 & mean & 2.5047 & -- & -- & -- & -- & -- \\ 
 & s.d. & 0.2872 & -- & -- & -- & -- & -- \\ \hline
2 & mean & 1.6112 & 0.2946 & -- & -- & -- & -- \\ 
 & s.d. & 0.2703 & 0.1986 & -- & -- & -- & -- \\ \hline
3 & mean & 2.7455 & -0.9144 & 2.3777 & -- & -- & -- \\ 
 & s.d. & 0.4856 & 0.3301 & 0.4465 & -- & -- & -- \\ \hline
$4^{\ast}$ & mean & 3.2248 & -1.0726 & 2.1506 & 1.0782 & -- & -- \\ 
 & s.d. & 0.6262 & 0.3792 & 0.5081 & 0.3492 & -- & -- \\ \hline
5 & mean & 3.2692 & -1.0873 & 2.1804 & 1.0957 & -0.0029 & -- \\ 
 & s.d. & 0.6508 & 0.3889 & 0.5250 & 0.4033 & 0.3008 & -- \\ \hline
6 & mean & 3.2998 & -1.0937 & 2.1975 & 1.1106 & -0.0045 & -0.0010 \\ 
 & s.d. & 0.6743 & 0.4020 & 0.5389 & 0.4199 & 0.3621 & 0.3081 \\ \hline
 \multicolumn{1}{l}{} & \multicolumn{1}{l}{} & \multicolumn{6}{c}{$n=300$} \\
  & & $\hat{\theta}_{1}$ & $\hat{\theta}_{2}$ & $\hat{\theta}_{3}$ & $\hat{\theta}_{4}$ & $\hat{\theta}_{5}$ & $\hat{\theta}_{6}$ \\ \hline
1 & mean & 2.5001 & -- & -- & -- & -- & -- \\ 
 & s.d. & 0.2878 & -- & -- & -- & -- & -- \\ \hline
2 & mean & 1.5930 & 0.2921 & -- & -- & -- & -- \\ 
 & s.d. & 0.2184 & 0.1615 & -- & -- & -- & -- \\ \hline
3 & mean & 2.6986 & -0.9016 & 2.3370 & -- & -- & -- \\ 
 & s.d. & 0.3753 & 0.2669 & 0.3494 & -- & -- & -- \\ \hline
$4^{\ast}$ & mean & 3.1478 & -1.0509 & 2.1000 & 1.0516 & -- & -- \\ 
 & s.d. & 0.4761 & 0.3026 & 0.3904 & 0.2718 & -- & -- \\ \hline
5 & mean & 3.1738 & -1.0590 & 2.1168 & 1.0628 & -0.0039 & -- \\ 
 & s.d. & 0.4864 & 0.3074 & 0.3980 & 0.3100 & 0.2321 & -- \\ \hline
6 & mean & 3.1914 & -1.0619 & 2.1249 & 1.0725 & -0.0090 & 0.0043 \\ 
 & s.d. & 0.5053 & 0.3177 & 0.4114 & 0.3209 & 0.2803 & 0.2353 \\ \hline
\end{tabular}
}
\label{esti3}
\end{table}


\section{Real data example} \label{real.simu}
The QBIC, BIC, and fAIC were also applied to the analysis of the meteorological data, which can be found at the Homepage of Japan Meteorological Agency (http://www.jma.go.jp/jma/indexe.html). 
The data was obtained during a period of 11 years from January 2000 to December 2010 at Yonagunijima, Japan. 
The data includes the monthly total precipitation $(P)_{t}$, monthly mean temperature $(T)_{t}$, monthly mean carbon dioxide $(CO_{2})_{t}$, monthly mean methane $(CH_{4})_{t}$, monthly mean carbon monoxide $(CO)_{t}$, and monthly mean ozone $(O_{3})_{t}$, where $t=-11,\ldots,-1,0,1,\ldots,120$.
The seasonal difference of precipitation is used for analysis, and $Y_{t}$, $t=1,\ldots,120$, are given by
\begin{align*}
Y_{t}=\left\{\begin{array}{cc}
1, & \text{if }(P)_{t}-(P)_{t-12}\geq0 \vspace{1mm}\\ 
0. & \text{if }(P)_{t}-(P)_{t-12}<0
\end{array}\right.
\end{align*}
In the parameter estimation and model selection, we use the linear logistic regression model
\begin{align*}
Y_{t}\sim B\bigg(1,\frac{\exp(X_{t-1}^{\prime}\theta)}{1+\exp(X_{t-1}^{\prime}\theta)}\bigg),
\end{align*}
and $X_{t}$ has the following elements in each candidate model:
\begin{align*}
&{\bf Model}\;{\bf1:}\;X_{t}=(X_{t,1},X_{t,2},X_{t,3},X_{t,4},X_{t,5}); \\
&{\bf Model}\;{\bf2:}\;X_{t}=(X_{t,1},X_{t,2},X_{t,3}X_{t,4});\;{\bf Model}\;{\bf3:}\;X_{t}=(X_{t,1},X_{t,2},X_{t,3},X_{t,5}); \\
&{\bf Model}\;{\bf4:}\;X_{t}=(X_{t,1},X_{t,2},X_{t,4},X_{t,5});\;{\bf Model}\;{\bf5:}\;X_{t}=(X_{t,1},X_{t,3},X_{t,4},X_{t,5}); \\
&\cdots;\;{\bf Model}\;{\bf26:}\;X_{t}=(X_{t,4},X_{t,5});\;\cdots;\;{\bf Model}\;{\bf31:}\;X_{t}=X_{t,5}.
\end{align*}
Here, $X_{t,1},X_{t,2},X_{t,3},X_{t,4}$, and $X_{t,5}$ denote the normalized $(T)_{t}$, $(CO_{2})_{t}$, $(CH_{4})_{t}$, $(CO)_{t}$, and $(O_{3})_{t}$.

The estimators and values of QBIC, BIC, and fAIC are shown in Table \ref{real.ex}.
By comparing the calculation results of QBIC and fAIC, Model 26, which consists of $X_{t,4}$ and $X_{t,5}$, is selected as the best model.
On the other hand, the calculation results of BIC imply that Model 31, which consists of $X_{t,5}$, is chosen.
Note that Model 26 contains Model 31.
These results mean that $(CO)_{t}$ and $(O_{3})_{t}$ are more significant than $(T)_{t}$, $(CO_{2})_{t}$, and $(CH_{4})_{t}$ for the seasonal difference of precipitation.

\begin{table}
\caption{\footnotesize The estimators and values of QBIC, BIC, and fAIC in each model (1--31 represent the models)}
\scalebox{0.8}[0.8]{
\begin{tabular}{ l c c c c c c c c c c } \hline
\multicolumn{1}{l}{} & & \multicolumn{5}{c}{Estimators} & & \multicolumn{3}{c}{Criteria} \\ \cline{3-7} \cline{9-11}
\multicolumn{1}{l}{} & & $(T)_{t}$ & $(CO_{2})_{t}$ & $(CH_{4})_{t}$ & $(CO)_{t}$ & $(O_{3})_{t}$ & & \multirow{2}{*}{QBIC} & \multirow{2}{*}{BIC} & \multirow{2}{*}{fAIC} \\
 & & $\hat{\theta}_{1}$ & $\hat{\theta}_{2}$ & $\hat{\theta}_{3}$ & $\hat{\theta}_{4}$ & $\hat{\theta}_{5}$ & & & & \\[1pt] \hline
1 & & -0.1328 & 0.0903 & -0.3482 & 0.5266 & -0.4888 & & 173.9913 & 186.2375 & 172.3000 \\[1pt] \hline 
2 & & -0.2615 & 0.1934 & -0.6884 & 0.2529 & -- & & 173.6681 & 182.6151 & 171.4651 \\[1pt] \hline 
3 & & -0.4314 & 0.0927 & -0.3403 & -- & -0.2720 & & 173.7480 & 182.4788 & 171.3288 \\[1pt] \hline 
4 & & -0.0210 & -0.0254 & -- & 0.5321 & -0.6850 & & 173.1102 & 181.8313 & 170.6813 \\[1pt] \hline 
5 & & -0.1137 & -- & -0.2385 & 0.5385 & -0.5630 & & 171.2799 & 181.5201 & 170.3702 \\[1pt] \hline 
6 & & -- & 0.0557 & -0.2633 & 0.6259 & -0.5396 & & 172.4866 & 181.5337 & 170.3838 \\[1pt] \hline 
7 & & -0.4119 & 0.1726 & -0.5915 & -- & -- & & 172.3932 & 178.1460 & 169.7835 \\[1pt] \hline 
8 & & -0.0167 & 0.0200 & -- & -0.0591 & -- & & 174.8113 & 180.6593 & 172.2968 \\[1pt] \hline 
9 & & -0.3163 & -0.0057 & -- & -- & -0.4567 & & 172.8990 & 178.1106 & 169.7482 \\[1pt] \hline 
10 & & -0.2325 & -- & -0.5170 & 0.1859 & -- & & 171.4395 & 178.6055 & 170.2430 \\[1pt] \hline 
11 & & -0.4111 & -- & -0.2185 & -- & -0.3449 & & 171.0992 & 177.8327 & 169.4702 \\[1pt] \hline 
12 & & -0.0154 & -- & -- & 0.5260 & -0.6796 & & 169.8614 & 177.0610 & 168.6985 \\[1pt] \hline 
13 & & -- & 0.1825 & -0.6038 & 0.4180 & -- & & 172.3837 & 178.1883 & 169.8259 \\[1pt] \hline 
14 & & -- & 0.0434 & 0.0045 & -- & -0.2272 & & 173.2589 & 179.2828 & 170.9203 \\[1pt] \hline 
15 & & -- & -0.0246 & -- & 0.5506 & -0.6846 & & 171.3291 & 177.0464 & 168.6839 \\[1pt] \hline 
16 & & -- & -- & -0.1973 & 0.6239 & -0.5848 & & 169.7155 & 176.7985 & 168.4360 \\[1pt] \hline 
17 & & 0.0357 & 0.0193 & -- & -- & -- & & 173.0224 & 175.8931 & 170.3181 \\[1pt] \hline 
18 & & -0.3488 & -- & -0.4575 & -- & -- & & 169.9879 & 173.9960 & 168.4211 \\[1pt] \hline 
19 & & -0.0210 & -- & -- & -0.0578 & -- & & 171.5117 & 175.8827 & 170.3077 \\[1pt] \hline 
20 & & -0.3142 & -- & -- & -- & -0.4557 & & 169.6164 & 173.3240 & 167.7491 \\[1pt] \hline 
21 & & -- & 0.1120 & -0.2194 & -- & -- & & 171.7523 & 174.8190 & 169.2441 \\[1pt] \hline 
22 & & -- & 0.0204 & -- & -0.0436 & -- & & 173.0122 & 175.8735 & 170.2985 \\[1pt] \hline 
23 & & -- & 0.0448 & -- & -- & -0.2237 & & 171.6428 & 174.4954 & 168.9204 \\[1pt] \hline 
24 & & -- & -- & -0.4511 & 0.3383 & -- & & 170.0575 & 174.1061 & 168.5311 \\[1pt] \hline 
25 & & -- & -- & 0.0548 & -- & -0.2628 & & 170.4487 & 174.5279 & 168.9530 \\[1pt] \hline 
26 & & -- & -- & -- & 0.5398 & -0.6794 & & {\bf 168.0681} & 172.2750 & {\bf 166.7000} \\[1pt] \hline 
27 & & 0.0302 & -- & -- & -- & -- & & 169.7206 & 171.1156 & 168.3281 \\[1pt] \hline 
28 & & -- & 0.0092 & -- & -- & -- & & 169.7456 & 171.1403 & 168.3528 \\[1pt] \hline 
29 & & -- & -- & -0.1669 & -- & -- & & 168.9109 & 170.3210 & 167.5335 \\[1pt] \hline 
30 & & -- & -- & -- & -0.0388 & -- & & 169.7026 & 171.0980 & 168.3105 \\[1pt] \hline 
31 & & -- & -- & -- & -- & -0.2165 & & 168.3457 & {\bf 169.7660} & 166.9785 \\[1pt] \hline
\end{tabular}
}
\label{real.ex}
\end{table}


\section{Proofs} \label{Proofs}

We will make use of the next three lemmas.
The proofs of Lemmas \ref{lem3}--\ref{lem2} are given in Supplementary Material.
Recall that $\psi_{j}$ is given by $\psi_{j}=\big(Y_{j}-F(X_{j})\big)X_{j}$ for all $j\in\mbbn$.

\begin{lem}
Assume that Assumptions \ref{Ass2} and \ref{Ass5} are satisfied and that $\sup_{j\in\mbbn}\|\psi_{j}\|_{2}<\infty$, then 
\begin{align*}
\sup_{n>0}\frac{1}{n}E\bigg[\sup_{1\leq i\leq n}\bigg|\sum_{j=1}^{i}\psi_{j}\bigg|^{2}\bigg]<\infty.
\end{align*}
\label{lem3}
\end{lem}

We write $\D_{n}=\frac{1}{\sqrt{n}}\p_{\theta}\mbbh_{n}(\tz)=\frac{1}{\sqrt{n}}\sumj\big(Y_{j}-\p b(X_{j}^{\prime}\tz)\big)X_{j}$.

\begin{lem}
Assume that Assumptions \ref{Ass1}--\ref{Ass3} and \ref{Ass7} are satisfied, then the following claims are established:
\begin{itemize}
\item[(i)] $\D_{n}=O_{p}(1)$.
\item[(ii)] $\displaystyle\sup_{\theta\in\Theta}\bigg|\frac{1}{n\sqrt{n}}\p_{\theta}^{3}\mbbh_{n}(\theta)\bigg|=o_{p}(1)$.
\end{itemize}
\label{lem1}
\end{lem}

We write $\mbbu_{n}(\tz)=\big\{u\in\mbbr^{p};\tz+\frac{u}{\sqrt{n}}\in\Theta\big\}$ and 
\begin{align*}
\mbbz_{n}(u)=\exp\Big\{\mbbh_{n}\Big(\tz+\frac{u}{\sqrt{n}}\Big)-\mbbh_{n}(\tz)\Big\}.
\end{align*}

\begin{lem}
If Assumptions \ref{Ass1}--\ref{Ass3} and \ref{Ass7} hold, then 
\begin{align*}
\int_{\mbbu_{n}(\tz)\cap\{|u|\geq M_{n}\}}\mbbz_{n}(u)du=o_{p}(1)
\end{align*}
for any $M_{n}\to\infty$.
\label{lem2}
\end{lem}


\subsection{Proof of Theorem \ref{qlfTh1}} \label{prfth1}
In what follows, we consider the zero-extended version of $\mbbz_{n}$ and use the same notation: 
\begin{align*}
\int_{\mbbr^{p}\backslash\mbbu_{n}(\tz)}\mbbz_{n}(u)du=0.
\end{align*}
By using the change of variable $\theta=\tz+\frac{u}{\sqrt{n}}$, the log marginal quasi-likelihood function becomes
\begin{align*}
&\log\bigg(\int_{\Theta}\exp\{\mbbh_{n}(\theta)\}\pi(\theta)d\theta\bigg) \\
&=\mbbh_{n}(\tz)-\frac{p}{2}\log n+\log\bigg\{\int_{\mbbu_{n}(\tz)}\mbbz_{n}(u)\pi\Big(\tz+\frac{u}{\sqrt{n}}\Big)du\bigg\} \\
&=\mbbh_{n}(\tz)-\frac{p}{2}\log n \\
&\qquad+\log\bigg\{\int_{\mbbu_{n}(\tz)}\mbbz_{n}(u)\bigg(\pi\Big(\tz+\frac{u}{\sqrt{n}}\Big)-\pi(\tz)\bigg)du+\pi(\tz)\int_{\mbbr^{p}}\mbbz_{n}(u)du\bigg\}.
\end{align*}

First we consider the asymptotic behavior of $\int_{\mbbu_{n}(\tz)}\mbbz_{n}(u)\big(\pi(\tz+\frac{u}{\sqrt{n}})-\pi(\tz)\big)du$. 
Because of the condition (ii) of Theorem \ref{qlfTh1}, Assumption \ref{Ass3} (i) and Lemma \ref{lem2}, we can take $M>0$ large enough so that
\begin{align*}
&\bigg|\int_{\mbbu_{n}(\tz)}\mbbz_{n}(u)\bigg(\pi\Big(\tz+\frac{u}{\sqrt{n}}\Big)-\pi(\tz)\bigg)du\bigg| \\
&\leq\int_{\mbbu_{n}(\tz)}\mbbz_{n}(u)\bigg|\pi\Big(\tz+\frac{u}{\sqrt{n}}\Big)-\pi(\tz)\bigg|du \\
&=\int_{\mbbu_{n}(\tz)\cap\{|u|<M\}}\mbbz_{n}(u)\bigg|\pi\Big(\tz+\frac{u}{\sqrt{n}}\Big)-\pi(\tz)\bigg|du \\
&\qquad+\int_{\mbbu_{n}(\tz)\cap\{|u|\geq M\}}\mbbz_{n}(u)\bigg|\pi\Big(\tz+\frac{u}{\sqrt{n}}\Big)-\pi(\tz)\bigg|du \\
&\leq(2M)^{p}\sup_{|u|<M}\bigg|\pi\Big(\tz+\frac{u}{\sqrt{n}}\Big)-\pi(\tz)\bigg|\sup_{|u|<M}\mbbz_{n}(u) \\
&\qquad+2\sup_{\theta\in\Theta}\pi(\theta)\int_{\mbbu_{n}(\tz)\cap\{|u|\geq M\}}\mbbz_{n}(u)du \\
&=o_{p}(1)\times\sup_{|u|<M}\bigg\{\exp\bigg(u^{\prime}\D_{n}-\frac{1}{2}u^{\prime}\bigg(\frac{1}{n}\sumj\p^{2} b(X_{j}^{\prime}\tilde{\theta}_{n})X_{j}X_{j}^{\prime}\bigg)u\bigg)\bigg\} \\
&\qquad+O_{p}(1)\times o_{p}(1) \\
&\leq o_{p}(1)\times\sup_{|u|<M}\bigg\{\exp\bigg(u^{\prime}\D_{n}-\frac{1}{2}u^{\prime}\bigg(\frac{1}{n}\sumj\underline{b}(X_{j})X_{j}X_{j}^{\prime}\bigg)u\bigg)\bigg\}+o_{p}(1),
\end{align*}
where $\tilde{\theta}_{n}=\tz+\xi \frac{u}{\sqrt{n}}$ for some $\xi$ satisfying $0<\xi<1$. 
Since $\frac{\p}{\p u}\big\{u^{\prime}\D_{n}-\frac{1}{2}u^{\prime}\big(\frac{1}{n}\sumj\underline{b}(X_{j})X_{j}X_{j}^{\prime}\big)u\big)\big\}=0$ if and only if $u=\big(\frac{1}{n}\sumj\underline{b}(X_{j})X_{j}X_{j}^{\prime}\big)^{-1}\D_{n}$, we have
\begin{align*}
u^{\prime}\D_{n}-\frac{1}{2}u^{\prime}\bigg(\frac{1}{n}\sumj\underline{b}(X_{j})X_{j}X_{j}^{\prime}\bigg)u&\leq\frac{1}{2}\D_{n}^{\prime}\bigg(\frac{1}{n}\sumj\underline{b}(X_{j})X_{j}X_{j}^{\prime}\bigg)^{-1}\D_{n}. 
\end{align*}
From Assumption \ref{Ass3} (ii) and Lemma \ref{lem1} (i), for any $\epsilon>0$ and for some $L>0$, 
\begin{align}
&\limsup_{n \to \infty}P\bigg[\sup_{|u|<M}\bigg\{\exp\bigg(u^{\prime}\D_{n}-\frac{1}{2}u^{\prime}\bigg(\frac{1}{n}\sumj\underline{b}(X_{j})X_{j}X_{j}^{\prime}\bigg)u\bigg)\bigg\}>L\bigg] \notag \\
&\leq\limsup_{n \to \infty}P\bigg[\exp\bigg\{\frac{1}{2}\D_{n}^{\prime}\bigg(\frac{1}{n}\sumj\underline{b}(X_{j})X_{j}X_{j}^{\prime}\bigg)^{-1}\D_{n}\bigg\}>L; \notag  \\
&\qquad\qquad\qquad\qquad\qquad\lambda_{\min}\bigg(\frac{1}{n}\sumj\underline{b}(X_{j})X_{j}X_{j}^{\prime}\bigg)<\lambda_{0}\bigg] \notag \\
&\qquad+\limsup_{n \to \infty}P\bigg[\exp\bigg\{\frac{1}{2}\D_{n}^{\prime}\bigg(\frac{1}{n}\sumj\underline{b}(X_{j})X_{j}X_{j}^{\prime}\bigg)^{-1}\D_{n}\bigg\}>L; \notag  \\
&\qquad\qquad\qquad\qquad\qquad\qquad\lambda_{\min}\bigg(\frac{1}{n}\sumj\underline{b}(X_{j})X_{j}X_{j}^{\prime}\bigg)\geq\lambda_{0}\bigg] \notag \\
&\leq\limsup_{n \to \infty}P\bigg[\lambda_{\min}\bigg(\frac{1}{n}\sumj\underline{b}(X_{j})X_{j}X_{j}^{\prime}\bigg)<\lambda_{0}\bigg] \notag  \\
&\qquad+\limsup_{n \to \infty}P\bigg[\exp\bigg\{\frac{1}{2\lambda_{0}}\D_{n}^{\prime}\D_{n}\bigg\}>L\bigg] \notag \\
&<\epsilon. \label{th_ine1}
\end{align}
Hence, $\sup_{|u|<M}\big\{\exp\big(u^{\prime}\D_{n}-\frac{1}{2}u^{\prime}\big(\frac{1}{n}\sumj\underline{b}(X_{j})X_{j}X_{j}^{\prime}\big)u\big)\big\}=O_{p}(1)$, and $\int_{\mbbu_{n}(\tz)}$ $\mbbz_{n}(u)\Big(\pi\big(\tz+\frac{u}{\sqrt{n}}\big)-\pi(\tz)\Big)$ converges to 0 in probability.

Next we will prove the equation $\int_{\mbbr^{p}}\mbbz_{n}(u)du=\int_{\mbbr^{p}}\exp(u^{\prime}\D_{n}-\frac{1}{2}u^{\prime}\Gam_{n}u)du+o_{p}(1)$. For each $K>0$,
\begin{align*}
&\bigg|\int_{\mbbr^{p}}\mbbz_{n}(u)du-\int_{\mbbr^{p}}\exp\Big(u^{\prime}\D_{n}-\frac{1}{2}u^{\prime}\Gam_{n}u\Big)du\bigg| \\
&\leq\int_{\mbbr^{p}}\bigg|\mbbz_{n}(u)-\exp\Big(u^{\prime}\D_{n}-\frac{1}{2}u^{\prime}\Gam_{n}u\Big)\bigg|du \\
&=\int_{|u|<K}\bigg|\mbbz_{n}(u)-\exp\Big(u^{\prime}\D_{n}-\frac{1}{2}u^{\prime}\Gam_{n}u\Big)\bigg|du \\
&\qquad+\int_{|u|\geq K}\bigg|\mbbz_{n}(u)-\exp\Big(u^{\prime}\D_{n}-\frac{1}{2}u^{\prime}\Gam_{n}u\Big)\bigg|du.
\end{align*}
Due to Assumption \ref{Ass3} and Lemma \ref{lem2}, we can take $K$ large enough so that
\begin{align*}
&\int_{|u|\geq K}\bigg|\mbbz_{n}(u)-\exp\Big(u^{\prime}\D_{n}-\frac{1}{2}u^{\prime}\Gam_{n}u\Big)\bigg|du \\
&\leq\int_{|u|\geq K}\mbbz_{n}(u)du \\
&\qquad+\int_{|u|\geq K}\exp\bigg\{u^{\prime}\D_{n}-\frac{1}{2}u^{\prime}\bigg(\frac{1}{n}\sumj\p^{2}b(X_{j}^{\prime}\tz)X_{j}X_{j}^{\prime}\bigg)u\bigg\}du \\
&\lesssim\int_{|u|\geq K}\exp\bigg\{u^{\prime}\D_{n}-\frac{1}{2}u^{\prime}\bigg(\frac{1}{n}\sumj\underline{b}(X_{j})X_{j}X_{j}^{\prime}\bigg)u\bigg\}du=o_{p}(1).
\end{align*}
In the same way as (\ref{th_ine1}), for the same $K$,
\begin{align*}
&\int_{|u|<K}\bigg|\mbbz_{n}(u)-\exp\Big(u^{\prime}\D_{n}-\frac{1}{2}u^{\prime}\Gam_{n}u\Big)\bigg|du \\
&\lesssim\sup_{|u|<K}\bigg|\exp\Big(u^{\prime}\D_{n}-\frac{1}{2}u^{\prime}\Gam_{n}u\Big) \\
&\qquad\qquad\times\bigg\{\exp\bigg(\frac{1}{6}\sum_{i,k,\ell=1}^{p}\frac{1}{n\sqrt{n}}\p_{\theta_{i}}\p_{\theta_{k}}\p_{\theta_{\ell}}\mbbh_{n}(\tilde{\theta}_{n})u_{i}u_{k}u_{\ell}\bigg)-1\bigg\}\bigg| \\
&=\sup_{|u|<K}\bigg\{\exp\bigg(u^{\prime}\D_{n}-\frac{1}{2}u^{\prime}\bigg(\frac{1}{n}\sumj\p^{2} b(X_{j}^{\prime}\tz)X_{j}X_{j}^{\prime}\bigg)u\bigg) \\
&\qquad\qquad\times\bigg|\exp\bigg(\frac{1}{6}\sum_{i,k,\ell=1}^{p}\frac{1}{n\sqrt{n}}\p_{\theta_{i}}\p_{\theta_{k}}\p_{\theta_{\ell}}\mbbh_{n}(\tilde{\theta}_{n})u_{i}u_{k}u_{\ell}\bigg)-1\bigg|\bigg\} \\
&\leq\sup_{|u|<K}\bigg\{\exp\bigg(u^{\prime}\D_{n}-\frac{1}{2}u^{\prime}\bigg(\frac{1}{n}\sumj\underline{b}(X_{j})X_{j}X_{j}^{\prime}\bigg)u\bigg) \\
&\qquad\qquad\times\bigg|\exp\bigg(\frac{1}{6}\sum_{i,k,\ell=1}^{p}\frac{1}{n\sqrt{n}}\p_{\theta_{i}}\p_{\theta_{k}}\p_{\theta_{\ell}}\mbbh_{n}(\tilde{\theta}_{n})u_{i}u_{k}u_{\ell}\bigg)-1\bigg|\bigg\} \\
&\leq\exp\bigg\{\frac{1}{2}\D_{n}^{\prime}\bigg(\frac{1}{n}\sumj\underline{b}(X_{j})X_{j}X_{j}^{\prime}\bigg)^{-1}\D_{n}\bigg\} \\
&\qquad\qquad\times\sup_{|u|<K}\bigg|\exp\bigg(\frac{1}{6}\sum_{i,k,\ell=1}^{p}\frac{1}{n\sqrt{n}}\p_{\theta_{i}}\p_{\theta_{k}}\p_{\theta_{\ell}}\mbbh_{n}(\tilde{\theta}_{n})u_{i}u_{k}u_{\ell}\bigg)-1\bigg| \\
&=O_{p}(1)\times o_{p}(1)=o_{p}(1).
\end{align*}
Therefore, we obtain the equation $\int_{\mbbr^{p}}\mbbz_{n}(u)du=\int_{\mbbr^{p}}\exp(u^{\prime}\D_{n}-\frac{1}{2}u^{\prime}\Gam_{n}u)du+o_{p}(1)$. 
Moreover,
\begin{align*}
&\int_{\mbbr^{p}}\exp\Big(u^{\prime}\D_{n}-\frac{1}{2}u^{\prime}\Gam_{n}u\Big)du \\
&=\exp\Big(\frac{1}{2}\|\Gam_{n}^{-\frac{1}{2}}\D_{n}\|^{2}\Big)\int_{\mbbr^{p}}\exp\Big(-\frac{1}{2}(u-\Gam_{n}^{-1}\D_{n})^{\prime}\Gam_{n}(u-\Gam_{n}^{-1}\D_{n})\Big)du \\
&=\exp\Big(\frac{1}{2}\|\Gam_{n}^{-\frac{1}{2}}\D_{n}\|^{2}\Big)(2\pi)^{\frac{p}{2}}\det(\Gam_{n})^{-\frac{1}{2}},
\end{align*}
hence $\log\big(\int_{\Theta}\exp\{\mbbh_{n}(\theta)\}\pi(\theta)d\theta\big)$ is given by
\begin{align*}
&\log\bigg(\int_{\Theta}\exp\{\mbbh_{n}(\theta)\}\pi(\theta)d\theta\bigg) \\
&=\mbbh_{n}(\tz)-\frac{p}{2}\log n+\log\Big\{\pi(\tz)\exp\Big(\frac{1}{2}\|\Gam_{n}^{-\frac{1}{2}}\D_{n}\|\Big)(2\pi)^{\frac{p}{2}}\det(\Gam_{n})^{-\frac{1}{2}}+o_{p}(1)\Big\} \\
&=\mbbh_{n}(\tz)-\frac{p}{2}\log n+\log\pi(\tz) \\
&\qquad\qquad+\frac{1}{2}\|\Gam_{0}^{-\frac{1}{2}}\D_{n}\|^{2}+\frac{p}{2}\log 2\pi+\log\det(\Gam_{n})^{-\frac{1}{2}}+o_{p}(1).
\end{align*}

Finally we replace $\tz$ with the QMLE $\tes$: 
\begin{align*}
\D_{n}&=\frac{1}{\sqrt{n}}\p_{\theta}\mbbh_{n}(\tz) \\
&=\big(\sqrt{n}(\tes-\tz)\big)^{\prime}\bigg(-\frac{1}{n}\p_{\theta}^{2}\mbbh_{n}(\check{\theta}_{n})\bigg),
\end{align*}
where $\check{\theta}_{n}=\tes+\eta_{1}(\tz-\tes)$ for some $\eta_{1}$ satisfying $0<\eta_{1}<1$. 
Because of Lemma \ref{lem1}, there exists a $\eta_{2}$ satisfying $0<\eta_{2}<1$ such that
\begin{align*}
-\frac{1}{n}\p_{\theta}^{2}\mbbh_{n}(\check{\theta}_{n})&=\Gam_{n}-\big(\sqrt{n}(\check{\theta}_{n}-\tz)\big)^{\prime}\bigg(\frac{1}{n\sqrt{n}}\p_{\theta}^{3}\mbbh_{n}\big(\tz+\eta_{2}(\check{\theta}_{n}-\tz)\big)\bigg) \\
&=\Gam_{n}-(1-\eta_{1})\big(\sqrt{n}(\tes-\tz)\big)^{\prime}\bigg(\frac{1}{n\sqrt{n}}\p_{\theta}^{3}\mbbh_{n}\big(\tz+\eta_{2}(\check{\theta}_{n}-\tz)\big)\bigg) \\
&=\Gam_{0}+o_{p}(1).
\end{align*}
Furthermore, we can show the equation $\frac{1}{n}\sumj b(X_{j}^{\prime}\tes)X_{j}X_{j}^{\prime}=\Gam_{n}+o_{p}(1)=\Gam_{0}+o_{p}(1)$ in a similar way, so we have
\begin{align*}
\mbbh_{n}(\tz)&=\mbbh_{n}(\tes)-\frac{1}{2}\big(\sqrt{n}(\tes-\tz)\big)^{\prime}\Gam_{0}\big(\sqrt{n}(\tes-\tz)\big)+o_{p}(1) \\
&=\mbbh_{n}(\tes)-\frac{1}{2}\big(\Gam_{0}^{-1}\D_{n}\big)^{\prime}\Gam_{0}\big(\Gam_{0}^{-1}\D_{n}\big)+o_{p}(1) \\
&=\mbbh_{n}(\tes)-\frac{1}{2}\|\Gam_{0}^{-\frac{1}{2}}\D_{n}\|^{2}+o_{p}(1).
\end{align*} 
Thus, the asymptotic behavior of the log marginal quasi-likelihood function is given by
\begin{align*}
\log\bigg(\int_{\Theta}\exp\{\mbbh_{n}(\theta)\}\pi(\theta)d\theta\bigg)&=\mbbh_{n}(\tes)-\frac{p}{2}\log n+\log\pi(\tes)+\frac{p}{2}\log 2\pi \\
&\qquad-\frac{1}{2}\log\det\bigg(\frac{1}{n}\sumj b(X_{j}^{\prime}\tes)X_{j}X_{j}^{\prime}\bigg)+o_{p}(1).
\end{align*}


\subsection{Proof of Theorem \ref{mod_consis}} \label{pr_mod_consis}
Recall that $\theta_{m,0}$ and $m_{0}$ are given by $\{\theta_{m,0}\}=\argmax_{\theta\in\Theta}\mbbh_{m,0}(\theta)$ and $\{m_{0}\}=\argmin_{m\in\mathcal{M}}\dim(\Theta_{m})$, respectively, where $\mathcal{M}=\argmax_{m\in\{1,\ldots,M\}}\mbbh_{m,0}$ $(\theta_{m,0})$. 
Fasen and Kimmig \cite{FaKim15} proved the model selection consistency of likelihood-based information criteria, which include AIC and BIC, for multivariate continuous-time ARMA processes. We basically follow their scenario for the proof of Theorem \ref{mod_consis}. \\
(i) $\Theta_{m_{0}}$ is nested in $\Theta_{m}$. Define the map $a:\Theta_{m_{0}}\to\Theta_{m}$ by $a(\theta)=A\theta+c$, where $A$ and $c$ satisfy the equation $\mbbh_{m_{0},n}(\theta)=\mbbh_{m,n}\big(a(\theta)\big)$ for any $\theta\in\Theta_{m_{0}}$.
Then, the equation $\mbbh_{m_{0},0}(\theta)=\mbbh_{m,0}\big(a(\theta)\big)$ is also satisfied for every $\theta\in\Theta_{m_{0}}$.
If $a(\theta_{m_{0},0})\neq \theta_{m,0}$, then we have $\mbbh_{m_{0},0}(\theta_{m_{0},0})=\mbbh_{m,0}\big(a(\theta_{m_{0},0})\big)<\mbbh_{m,0}(\theta_{m,0})$.
From this inequality and assumption of the optimal model, we have $a(\theta_{m_{0},0})=\theta_{m,0}$.

By the Taylor expansion of $\mbbh_{m,n}$
\begin{align*}
\mbbh_{m_{0},n}(\hat{\theta}_{m_{0},n})&=\mbbh_{m,n}\big(a(\hat{\theta}_{m_{0},n})\big) \\
&=\mbbh_{m,n}(\hat{\theta}_{m,n})-\frac{1}{2}\big\{\sqrt{n}\big(\hat{\theta}_{m,n}-a(\hat{\theta}_{m_{0},n})\big)\big\}^{\prime} \\
&\qquad\qquad\quad\times\bigg(\frac{1}{n}\sumj\p^{2}b_{m}(X_{j}^{\prime}\tilde{\theta}_{n})X_{j}X_{j}^{\prime}\bigg)\big\{\sqrt{n}\big(\hat{\theta}_{m,n}-a(\hat{\theta}_{m_{0},n})\big)\big\},
\end{align*}
where $\tilde{\theta}_{n}=\hat{\theta}_{m,n}+\xi\big(a(\hat{\theta}_{m_{0},n})-\hat{\theta}_{m,n}\big)$ for some $\xi$ satisfying $0<\xi<1$ and $\tilde{\theta}_{n}\cip\theta_{m,0}$ as $n\to\infty$.
Therefore, the difference between $\qbic^{(m_{0})}$ and $\qbic^{(m)}$ is given by
\begin{align*}
&\qbic^{(m_{0})}-\qbic^{(m)} \\
&=\big\{\sqrt{n}\big(\hat{\theta}_{m,n}-a(\hat{\theta}_{m_{0},n})\big)\big\}^{\prime}\bigg(\frac{1}{n}\sumj\p^{2}b_{m}(X_{j}^{\prime}\tilde{\theta}_{n})X_{j}X_{j}^{\prime}\bigg)\big\{\sqrt{n}\big(\hat{\theta}_{m,n}-a(\hat{\theta}_{m_{0},n})\big)\big\} \\
&\qquad+\log\det\big(-\p_{\theta}^{2}\mbbh_{m_{0},n}(\hat{\theta}_{m_{0},n})\big)-\log\det\big(-\p_{\theta}^{2}\mbbh_{m,n}(\hat{\theta}_{m,n})\big).
\end{align*} 

We consider the behavior of the $\hat{\theta}_{m,n}-a(\hat{\theta}_{m_{0},n})$. 
Because of the chain rule, we have $\p_{\theta}\mbbh_{m_{0},n}(\theta_{m_{0},0})=A^{\prime}\p_{\theta}\mbbh_{m,n}(\theta_{m,0})$ and $\p_{\theta}^{2}\mbbh_{m_{0},n}(\theta)=A^{\prime}\p_{\theta}^{2}\mbbh_{m,n}\big(a(\theta)\big)A$.
Moreover,
\begin{align*}
a(\hat{\theta}_{m_{0},n})-\theta_{m,0}&=A(\hat{\theta}_{m_{0},n}-\theta_{m_{0},0})
\end{align*}
and
\begin{align*}
&\sqrt{n}(\hat{\theta}_{m_{0},n}-\theta_{m_{0},0}) \\
&=\bigg(-\frac{1}{n}\p_{\theta}^{2}\mbbh_{m_{0},n}(\check{\theta}_{n})\bigg)^{-1}\bigg(\frac{1}{\sqrt{n}}\p_{\theta}\mbbh_{m_{0},n}(\theta_{m_{0},0})\bigg) \\
&=\bigg\{A^{\prime}\bigg(-\frac{1}{n}\p_{\theta}^{2}\mbbh_{m,n}\big(a(\check{\theta}_{n})\big)\bigg)A\bigg\}^{-1}A^{\prime}\bigg(\frac{1}{\sqrt{n}}\p_{\theta}\mbbh_{m,n}(\theta_{m,0})\bigg) \\
&=\bigg\{A^{\prime}\bigg(\frac{1}{n}\sumj\p^{2}b_{m}\big(a(\check{\theta}_{n})\big)X_{j}X_{j}^{\prime}\bigg)A\bigg\}^{-1}A^{\prime}\bigg(\frac{1}{\sqrt{n}}\p_{\theta}\mbbh_{m,n}(\theta_{m,0})\bigg) \\
&=\Big(A^{\prime}\Gam_{m,0}A\Big)^{-1}A^{\prime}\bigg(\frac{1}{\sqrt{n}}\p_{\theta}\mbbh_{m,n}(\theta_{m,0})\bigg)+o_{p}(1),
\end{align*}
where $\check{\theta}_{n}=\hat{\theta}_{m_{0},n}+\eta\big(\theta_{m_{0},0}-\hat{\theta}_{m_{0},n}\big)$ for some $\eta$ satisfying $0<\eta<1$ and $a(\check{\theta}_{n})\cip a(\theta_{m_{0},0})=\theta_{m,0}$ as $n\to\infty$. 
These equalities and Theorem \ref{AN} give
\begin{align*}
&\sqrt{n}(\hat{\theta}_{m,n}-a(\hat{\theta}_{m_{0},n})) \\
&=\sqrt{n}(\hat{\theta}_{m,n}-\theta_{m,0})-A\sqrt{n}(\hat{\theta}_{m_{0},n}-\theta_{m_{0},0}) \\
&\cil\Big\{\Gam_{m,0}^{-1}-A\Big(A^{\prime}\Gam_{m,0}A\Big)^{-1}A^{\prime}\Big\}N_{p_{m}}(0,\Sig_{0}) \\
&=N_{p_{m}}\bigg(0,\Big\{\Gam_{m,0}^{-1}-A\Big(A^{\prime}\Gam_{m,0}A\Big)^{-1}A^{\prime}\Big\}\Sig_{0}\Big\{\Gam_{m,0}^{-1}-A\Big(A^{\prime}\Gam_{m,0}A\Big)^{-1}A^{\prime}\Big\}\bigg)\sim{\bf N}.
\end{align*}
Thus,
\begin{align*}
& P[\qbic^{(m_{0})}-\qbic^{(m)}<0] \\
&=P\bigg[{\bf N}^{\prime}\bigg(\frac{1}{n}\sumj\p^{2}b_{m}(X_{j}^{\prime}\tilde{\theta}_{n})X_{j}X_{j}^{\prime}\bigg){\bf N}+\log\det\bigg(\frac{1}{n}\sumj\p^{2}b_{m_{0}}\big(X_{j}^{\prime}\hat{\theta}_{m_{0},n}\big)X_{j}X_{j}^{\prime}\bigg) \\
&\qquad-\log\det\bigg(\frac{1}{n}\sumj\p^{2}b_{m}(X_{j}^{\prime}\hat{\theta}_{m,n})X_{j}X_{j}^{\prime}\bigg)<(p_{m}-p_{m_{0}})\log n\bigg] \\
&\to P\Big[{\bf N}^{\prime}\Gam_{m,0}{\bf N}+\log\det\big(\Gam_{m_{0},0}\big)-\log\det\big(\Gam_{m,0}\big)<\infty\Big]
\end{align*}
as $n\to\infty$. 
From Imhof \cite[(1.1)]{Imh61}, ${\bf N}^{\prime}\Gam_{m,0}{\bf N}=\sum_{j=1}^{p_{m}}\lambda_{j}\chi_{j}^{2}$ in distribution, where ($\chi_{j}^{2}$) is a sequence of independent $\chi^{2}$ random variables with one degree of freedom and $\lambda_{j}$ are the eigenvalues of $\Gam_{m,0}^{\frac{1}{2}}\big\{\Gam_{m,0}^{-1}-A\big(A^{\prime}\Gam_{m,0}A\big)^{-1}A^{\prime}\big\}\Sig_{0}\big\{\Gam_{m,0}^{-1}-A\big(A^{\prime}\Gam_{m,0}A\big)^{-1}A^{\prime}\big\}\Gam_{m,0}^{\frac{1}{2}}$. 
Furthermore, $\log\det\big(\Gam_{m_{0},0}\big)=O(1)$ and $\log\det\big(\Gam_{m,0}\big)$ $=O(1)$.
Hence, 
\begin{align*}
& P\big[{\bf N}^{\prime}\Gam_{m,0}{\bf N}+\log\det\big(\Gam_{m_{0},0}\big)-\log\det\big(\Gam_{m,0}\big)<\infty\big] \\
&\geq P\bigg[\max_{j\in\{1,\ldots,p_{m}\}}\lambda_{j}\sum_{j=1}^{p_{m}}\chi_{j}^{2}<\infty\bigg]=1.
\end{align*}
(ii) $\mbbh_{m,0}(\theta)\neq\mbbh_{m_{0},0}(\theta_{m_{0},0})$ for every $\theta\in\Theta_{m}$. 
Because of Lemma \ref{lem1} (i) and the consistency of $\hat{\theta}_{m_{0},n}$ and $\hat{\theta}_{m,n}$, we have
\begin{align*}
\frac{1}{n}\mbbh_{m_{0},n}(\hat{\theta}_{m_{0},n})&=\frac{1}{n}\mbbh_{m_{0},n}(\theta_{m_{0},0})+o_{p}(1)=\mbbh_{m_{0},0}(\theta_{m_{0},0})+o_{p}(1)
\end{align*}
and
\begin{align*}
\frac{1}{n}\mbbh_{m,n}(\hat{\theta}_{m,n})&=\frac{1}{n}\mbbh_{m,n}(\theta_{m,0})+o_{p}(1)=\mbbh_{m,0}(\theta_{m,0})+o_{p}(1).
\end{align*}
Since $\mbbh_{m_{0},0}(\theta_{m_{0},0})$ is lager than $\mbbh_{m,0}(\theta_{m,0})$, we obtain
\begin{align*}
& P[\qbic^{(m_{0})}-\qbic^{(m)}<0] \\
&=P\bigg[-2\mbbh_{m_{0},n}(\hat{\theta}_{m_{0},n})+2\mbbh_{m,n}(\hat{\theta}_{m,n})+\log\det\bigg(\frac{1}{n}\sumj\p^{2}b_{m}\big(X_{j}^{\prime}\hat{\theta}_{m_{0},n}\big)X_{j}X_{j}^{\prime}\bigg) \\
&\qquad\qquad-\log\det\bigg(\frac{1}{n}\sumj\p^{2}b_{m}(X_{j}^{\prime}\hat{\theta}_{m,n})X_{j}X_{j}^{\prime}\bigg)<(p_{m}-p_{m_{0}})\log n\bigg] \\
&=P\bigg[\frac{-2}{n}\big(\mbbh_{m_{0},n}(\hat{\theta}_{m_{0},n})-\mbbh_{m,n}(\hat{\theta}_{m,n})\big) \\
&\qquad\qquad+\frac{1}{n}\log\det\bigg(\frac{1}{n}\sumj\p^{2}b_{m}\big(X_{j}^{\prime}\hat{\theta}_{m_{0},n}\big)X_{j}X_{j}^{\prime}\bigg) \\
&\qquad\qquad-\frac{1}{n}\log\det\bigg(\frac{1}{n}\sumj\p^{2}b_{m}(X_{j}^{\prime}\hat{\theta}_{m,n})X_{j}X_{j}^{\prime}\bigg)<(p_{m}-p_{m_{0}})\frac{\log n}{n}\bigg] \\
&\to P\Big[-2\big(\mbbh_{m_{0},0}(\theta_{m_{0},0})-\mbbh_{m,0}(\theta_{m,0})\big)<0\Big]=1
\end{align*}
as $n\to\infty$.


\section*{Acknowledgements}
The author wishes to thank the associate editor and the two anonymous referees for careful reading and valuable comments which helped to greatly improve the paper.
The author also thanks Professor H. Masuda for helpful comments and stimulating discussion.
This work was partly supported by JST, CREST.




\section{Supplementary Material}


\subsection{Proof of Lemma \ref{lem3}} \label{SM1}
Lemma \ref{lem3} follows from a direct application of Yoshida \cite[Lemma 4]{Yos11}.


\subsection{Proof of Lemma \ref{lem1}} \label{SM2}
\textit{Proof of (i).} 
\begin{align}
\sup_{j\in\mbbn}E[|\psi_{j}|^{2}]&=\sup_{j\in\mbbn}E\big[\big|\big(Y_{j}-F(X_{j})\big)X_{j}\big|^{2}\big] \notag \\
&\leq\sup_{j\in\mbbn}E\big[\big|Y_{j}-F(X_{j})\big|^{2}\big|X_{j}\big|^{2}\big] \notag \\
&\lesssim\sup_{j\in\mbbn}E\big[\big(|Y_{j}|^{2}+\big|F(X_{j})\big|^{2}\big)\big|X_{j}\big|^{2}\big] \notag \\
&\leq\sup_{j\in\mbbn}E\big[\big(|Y_{j}|^{2}+E[|Y_{j}|^{2}|\mcf_{j-1}\vee\sigma(X_{j})]\big)\big|X_{j}\big|^{2}\big] \notag \\
&\lesssim\sup_{j\in\mbbn}E\big[(1+|X_{j}|^{C^{\prime}})\big|X_{j}\big|^{2}\big]<\infty. \label{lem1_ineq1}
\end{align}
Because of this inequality, we can apply Lemma \ref{lem3} to obtain 
\begin{align*} 
\sup_{n>0}E\bigg[\bigg|\frac{1}{\sqrt{n}}\sumj\psi_{j}\bigg|^{2}\bigg]&\leq\sup_{n>0}\frac{1}{n}E\bigg[\sup_{1\leq i\leq n}\bigg|\sum_{j=1}^{i}\psi_{j}\bigg|^{2}\bigg]<\infty.
\end{align*}
Therefore, $\frac{1}{\sqrt{n}}\sumj\psi_{j}=O_{p}(1)$, and $\D_{n}$ satisfies the equality
\begin{align*}
\D_{n}=\frac{1}{\sqrt{n}}\sumj\psi_{j}+\frac{1}{\sqrt{n}}\sumj\big(F(X_{j})-\p b(X_{j}^{\prime}\tz)\big)X_{j}=O_{p}(1).
\end{align*}
\textit{Proof of (ii).} For some $C>0$,
\begin{align*}
\sup_{\theta\in\Theta}\bigg|\frac{1}{n\sqrt{n}}\p_{\theta}^{3}\mbbh_{n}(\theta)\bigg|&\leq\frac{1}{n\sqrt{n}}\sup_{\theta\in\Theta}\bigg(\sum_{i,k,\ell=1}^{p}\big|\p_{\theta_{i}}\p_{\theta_{k}}\p_{\theta_{\ell}}\mbbh_{n}(\theta)\big|^{2}\bigg)^{\frac{1}{2}} \\
&\lesssim\frac{1}{n\sqrt{n}}\sup_{\theta\in\Theta}\sum_{i,k,\ell=1}^{p}\bigg|\sumj\p^{3}b(X_{j}^{\prime}\theta)X_{j,i}X_{j,k}X_{j,\ell}\bigg| \\
&\leq\frac{1}{n\sqrt{n}}\sup_{\theta\in\Theta}\sum_{i,k,\ell=1}^{p}\sumj\big|\p^{3}b(X_{j}^{\prime}\theta)\big|\big|X_{j,i}X_{j,k}X_{j,\ell}\big| \\
&\leq\sum_{i,k,\ell=1}^{p}\frac{1}{n\sqrt{n}}\sumj\sup_{\theta\in\Theta}\big|\p^{3}b(X_{j}^{\prime}\theta)\big|\big|X_{j}\big|^{3} \\
&\lesssim\sum_{i,k,\ell=1}^{p}\frac{1}{n\sqrt{n}}\sumj\big(1+|X_{j}|^{C}\big)|X_{j}|^{3} \\
&=\sum_{i,k,\ell=1}^{p}\frac{1}{n\sqrt{n}}\sumj O_{p}(1)=o_{p}(1).
\end{align*}


\subsection{Proof of Lemma \ref{lem2}} \label{SM3}
We have that
\begin{align*}
&\int_{\mbbu_{n}(\tz)\cap\{|u|\geq M_{n}\}}\mbbz_{n}(u)du \\
&=\int_{\mbbu_{n}(\tz)\cap\{|u|\geq M_{n}\}}\exp\bigg\{u^{\prime}\D_{n}+\frac{1}{2n}u^{\prime}\p_{\theta}^{2}\mbbh_{n}(\tilde{\theta}_{n})u\bigg\}du \\
&=\int_{\mbbu_{n}(\tz)\cap\{|u|\geq M_{n}\}}\exp\bigg\{u^{\prime}\D_{n}-\frac{1}{2}u^{\prime}\bigg(\frac{1}{n}\sumj\p^{2}b(X_{j}^{\prime}\tilde{\theta}_{n})X_{j}X_{j}^{\prime}\bigg)u\bigg\}du,
\end{align*}
where $\tilde{\theta}_{n}=\tz+\xi(\tz+\frac{u}{\sqrt{n}}-\tz)=\tz+\xi \frac{u}{\sqrt{n}}$ for some $\xi$ satisfying $0<\xi<1$. 
From Assumption \ref{Ass3} (i), there exists a function $\underline{b}$ such that
\begin{align*}
& \int_{\mbbu_{n}(\tz)\cap\{|u|\geq M_{n}\}}\exp\bigg\{u^{\prime}\D_{n}-\frac{1}{2}u^{\prime}\bigg(\frac{1}{n}\sumj\p^{2}b(X_{j}^{\prime}\tilde{\theta}_{n})X_{j}X_{j}^{\prime}\bigg)u\bigg\}du \\
&\leq\int_{\mbbu_{n}(\tz)\cap\{|u|\geq M_{n}\}}\exp\bigg\{u^{\prime}\D_{n}-\frac{1}{2}u^{\prime}\bigg(\frac{1}{n}\sumj\underline{b}(X_{j})X_{j}X_{j}^{\prime}\bigg)u\bigg\}du \\
&\leq\int_{|u|\geq M_{n}}\exp\bigg\{u^{\prime}\D_{n}-\frac{1}{2}u^{\prime}\bigg(\frac{1}{n}\sumj\underline{b}(X_{j})X_{j}X_{j}^{\prime}\bigg)u\bigg\}du.
\end{align*}  
Fix any $\epsilon>0$. For $\lambda_{0}>0$ given in Assumption \ref{Ass3} (ii),
\begin{align}
& P\bigg[\int_{\mbbu_{n}(\tz)\cap\{|u|\geq M_{n}\}}\mbbz_{n}(u)du>\epsilon\bigg] \notag\\
&\leq P\bigg[\int_{|u|\geq M_{n}}\exp\bigg\{u^{\prime}\D_{n}-\frac{1}{2}u^{\prime}\bigg(\frac{1}{n}\sumj\underline{b}(X_{j})X_{j}X_{j}^{\prime}\bigg)u\bigg\}du>\epsilon\bigg] \notag\\
&=P\bigg[\int_{|u|\geq M_{n}}\exp\bigg\{u^{\prime}\D_{n}-\frac{1}{2}u^{\prime}\bigg(\frac{1}{n}\sumj\underline{b}(X_{j})X_{j}X_{j}^{\prime}\bigg)u\bigg\}du>\epsilon; \notag\\
&\qquad\qquad\qquad\lambda_{\min}\bigg(\frac{1}{n}\sumj\underline{b}(X_{j})X_{j}X_{j}^{\prime}\bigg)<\lambda_{0}\bigg] \notag\\
&\qquad+P\bigg[\int_{|u|\geq M_{n}}\exp\bigg\{u^{\prime}\D_{n}-\frac{1}{2}u^{\prime}\bigg(\frac{1}{n}\sumj\underline{b}(X_{j})X_{j}X_{j}^{\prime}\bigg)u\bigg\}du>\epsilon; \notag\\
&\qquad\qquad\qquad\qquad\lambda_{\min}\bigg(\frac{1}{n}\sumj\underline{b}(X_{j})X_{j}X_{j}^{\prime}\bigg)\geq\lambda_{0}\bigg] \notag\\
&\leq P\bigg[\lambda_{\min}\bigg(\frac{1}{n}\sumj\underline{b}(X_{j})X_{j}X_{j}^{\prime}\bigg)<\lambda_{0}\bigg] \notag\\
&\qquad+P\bigg[\int_{|u|\geq M_{n}}\exp\bigg\{u^{\prime}\D_{n}-\frac{1}{2}\lambda_{0}u^{\prime}u\bigg\}du>\epsilon\bigg]. \label{lem2_ineq1}
\end{align}
There exists a constant $K>0$ such that
\begin{align*}
& P\bigg[\int_{|u|\geq M_{n}}\exp\bigg\{u^{\prime}\D_{n}-\frac{1}{2}\lambda_{0}u^{\prime}u\bigg\}du>\epsilon\bigg] \\
&=P\bigg[\int_{|u|\geq M_{n}}\exp\bigg\{u^{\prime}\D_{n}-\frac{1}{2}\lambda_{0}u^{\prime}u\bigg\}du>\epsilon;|\D_{n}|>K\bigg] \notag \\
&\qquad\qquad+P\bigg[\int_{|u|\geq M_{n}}\exp\bigg\{u^{\prime}\D_{n}-\frac{1}{2}\lambda_{0}u^{\prime}u\bigg\}du>\epsilon;|\D_{n}|\leq K\bigg] \\
&\leq P\Big[|\D_{n}|>K\Big]+P\bigg[\exp\bigg(\frac{\D_{n}^{\prime}\D_{n}}{2\lambda_{0}}\bigg) \\
&\qquad\qquad\times\int_{|u|\geq M_{n}}\exp\bigg\{-\frac{\lambda_{0}}{2}(u-\lambda_{0}^{-1}\D_{n})^{\prime}(u-\lambda_{0}^{-1}\D_{n})\bigg\}du>\epsilon;|\D_{n}|\leq K\bigg] \\
&=P\Big[|\D_{n}|>K\Big]+P\bigg[\exp\bigg(\frac{\D_{n}^{\prime}\D_{n}}{2\lambda_{0}}\bigg) \\
&\qquad\qquad\times\int_{|t+\lambda_{0}^{-1}\D_{n}|\geq M_{n}}\exp\bigg\{-\frac{\lambda_{0}}{2}t^{\prime}t\bigg\}dt>\epsilon;|\D_{n}|\leq K\bigg] \\
&\leq P\Big[|\D_{n}|>K\Big]+P\bigg[\exp\bigg(\frac{K^{2}}{2\lambda_{0}}\bigg)\int_{|t|\geq M_{n}-\lambda_{0}^{-1}K}\exp\bigg\{-\frac{\lambda_{0}}{2}t^{\prime}t\bigg\}dt>\epsilon\bigg].
\end{align*}
Because of Lemma \ref{lem1} (i), for some $N$, 
\begin{align}
&P\Big[|\D_{n}|>K\Big] \notag\\
&\qquad+P\bigg[\exp\bigg(\frac{K^{2}}{2\lambda_{0}}\bigg)\int_{|t|\geq M_{n}-\lambda_{0}^{-1}K}\exp\bigg\{-\frac{\lambda_{0}}{2}t^{\prime}t\bigg\}dt>\epsilon\bigg]<\frac{\epsilon}{2} \label{lem2_ineq2}
\end{align}
for every $n\geq N$.
Due to Assumption \ref{Ass3} (ii), (\ref{lem2_ineq1}), and (\ref{lem2_ineq2}),
\begin{align*}
P\bigg[\int_{\mbbu_{n}(\tz)\cap\{|u|\geq M_{n}\}}\mbbz_{n}(u)du>\epsilon\bigg]<\epsilon
\end{align*}
for all $n\geq N$. Thus, $\displaystyle\int_{\mbbu_{n}(\tz)\cap\{|u|\geq M_{n}\}}\mbbz_{n}(u)du$ converges to 0 in probability.


\subsection{Proof of Theorem \ref{consis}} \label{pr_consis}
We apply the argmax theorem (van der Vaart \cite[Theorem 5.56, Corollary 5.58]{Van}) for $\frac{1}{n}\mbbh_{n}(\theta)$.
Under Assumptions \ref{Ass1}--\ref{Ass3}, $\tes$ and $\tz$ are given by $\{\tes\}=\argmax_{\theta\in\Theta}\frac{1}{n}\mbbh_{n}(\theta)$ and $\{\tz\}=\argmax_{\theta\in\Theta}\mbbh_{0}(\theta)$, respectively.
Hence, it is enough to show that $\frac{1}{n}\mbbh_{n}$ converges to $\mbbh_{0}$ uniformly in $\theta\in\Theta$.
Since (\ref{unicon1}) and Lemma \ref{lem3} are satisfied, and $\Theta$ is a bounded convex domain, we have
\begin{align*}
&\sup_{\theta\in\Theta}\left|\frac{1}{n}\mbbh_{n}(\theta)-\mbbh_{0}(\theta)\right| \\
&=\sup_{\theta\in\Theta}\bigg|\frac{1}{n}\sumj\psi_{j}^{\prime}\theta \\
&\qquad+\frac{1}{n}\sumj\big\{\big(F(X_{j})X_{j}^{\prime}\theta-b(X_{j}^{\prime}\theta)\big)\big\}-\int\big(F(x)x^{\prime}\theta-b(x^{\prime}\theta)\big)\nu(dx)\bigg| \\
&\leq\left|\frac{1}{n}\sumj\psi_{j}\right|\sup_{\theta\in\Theta}|\theta| \\
&\qquad+\sup_{\theta\in\Theta}\left|\frac{1}{n}\sumj\big\{\big(F(X_{j})X_{j}^{\prime}\theta-b(X_{j}^{\prime}\theta)\big)\big\}-\int\big(F(x)x^{\prime}\theta-b(x^{\prime}\theta)\big)\nu(dx)\right| \\
&\cip0.
\end{align*}


\subsection{Proof of Theorem \ref{AN}} \label{pr_AN}
It was shown in Section \ref{prfth1} that
\begin{align*}
\sqrt{n}(\tes-\tz)=\Gam_{0}^{-1}\D_{n}+o_{p}(1).
\end{align*}
In view of Herrndorf \cite[Theorem, Corollary 1]{Her84}, $\D_{n}$ converges to the normal distribution $N(0,\Sig_{0})$ in law if we show the following four conditions:
\begin{itemize}
\item[(i)] $E\big[\big(Y_{j}-\p b(X_{j}^{\prime}\tz)\big)X_{j}\big]=0$, $E\big[\big(Y_{j}-\p b(X_{j}^{\prime}\tz)\big)^{2}X_{j}^{\prime}X_{j}\big]<\infty$ for all $j\in\mbbn$.
\item[(ii)] $\frac{1}{n}E\Big[\Big\{\sumj\big(Y_{j}-\p b(X_{j}^{\prime}\tz)\big)X_{j}\Big\}\Big\{\sumj\big(Y_{j}-\p b(X_{j}^{\prime}\tz)\big)X_{j}\Big\}^{\prime}\Big]\to\Sig_{0}$ as $n\to\infty$.
\item[(iii)] $\sum_{k\in\mbbn}\alpha(k)^{\frac{1}{3}}<\infty$.
\item[(iv)] $\lim\sup_{n\to\infty}E\big[\big|\big(Y_{n}-\p b(X_{n}^{\prime}\tz)\big)X_{n}\big|^{3}\big]<\infty$.
\end{itemize}
(ii) is ensured by Assumption \ref{Ass6} (ii). \\
\textit{Proof of (i).} Because of Assumptions \ref{Ass2}, \ref{Ass6} (i), and the definition of $\tz$, we have
\begin{align*}
E\big[\big(Y_{j}-\p b(X_{j}^{\prime}\tz)\big)X_{j}\big]&=E[\big(Y_{j}-F(X_{j})\big)X_{j}+\big(F(X_{j})-\p b(X_{j}^{\prime}\tz)\big)X_{j}] \\
&=0+\int\big(F(x)-\p b(x^{\prime}\tz)\big)x\nu(dx) \\
&=0
\end{align*}
for any $j\in\mbbn$. Furthermore, from Assumption \ref{Ass1},
\begin{align}
\sup_{j\in\mbbn}E\big[\big(Y_{j}-\p b(X_{j}^{\prime}\tz)\big)^{2}X_{j}^{\prime}X_{j}\big]&\lesssim\sup_{j\in\mbbn}E\Big[\Big(|Y_{j}|^{2}+|\p b(X_{j}^{\prime}\tz)|^{2}\Big)|X_{j}|^{2}\Big] \notag \\
&\lesssim \sup_{j\in\mbbn}E\Big[\Big((1+|X_{j}|^{C^{\prime}})+(1+|X_{j}|^{C})^{2}\Big)|X_{j}|^{2}\Big] \notag \\
&<\infty. \label{pr_AN_ineq1}
\end{align} 
\textit{Proof of (iii).} Assumption \ref{Ass5} gives
\begin{align*}
\sum_{k\in\mbbn}\alpha(k)^{\frac{1}{3}}&\leq c^{-\frac{1}{3}}\sum_{k\in\mbbn}e^{-\frac{1}{3}ck} \\
&=c^{-\frac{1}{3}}\frac{e^{-\frac{1}{3}c}}{1-e^{-\frac{1}{3}c}}<\infty.
\end{align*}
\textit{Proof of (iv).} In a similar way as (\ref{pr_AN_ineq1}), we can show that
\begin{align*}
\sup_{j\in\mbbn}E\Big[\big|\big(Y_{j}-\p b(X_{j}^{\prime}\tz)\big)X_{j}\big|^{3}\Big]&<\infty
\end{align*}
Hence, (iv) is satisfied.

\end{document}